\documentclass{article}

\usepackage{arxiv}

\usepackage[utf8]{inputenc} % allow utf-8 input
\usepackage[T1]{fontenc}    % use 8-bit T1 fonts
\usepackage{hyperref}       % hyperlinks
\usepackage{url}            % simple URL typesetting
\usepackage{booktabs}       % professional-quality tables
\usepackage{amsfonts}       % blackboard math symbols
\usepackage{nicefrac}       % compact symbols for 1/2, etc.
\usepackage{microtype}      % microtypography
\usepackage{lipsum}
\usepackage{graphicx}

\title{A Multihorizon Approach for the Reliability Oriented Network Restructuring Problem Considering Learning Effects, Construction Time, and Cables Maintenance Costs}

\author{
  Chiara Bordin\thanks{Corresponding Author: Chiara Bordin; email: chiara.bordin@uit.no} \\
  Department of Computer Science\\
  UiT, The Arctic University of Norway \\
   \And
 Sambeet Mishra  \\
   Department of Electrical Power Engineering and Mechatronics\\
  Tallinn University of Technology \\
   \And
 Ivo Palu  \\
   Department of Electrical Power Engineering and Mechatronics\\
  Tallinn University of Technology \\
}

\begin{document}
\maketitle

\begin{abstract}
This paper presents a techno-economic optimisation tool to study how the power system expansion decisions can be taken in a more economical and efficient way, by minimising the consequent costs of network reinforcement and reconfiguration. Analyses are performed to investigate how the network reinforcement and reconfiguration should be planned, within a time horizon of several years, by continuously keeping the network feasibility and ability to satisfy the load.
The main contribution of this study is the inclusion of key features within the mathematical model to enhance the investment decision making process. A representative maintenance cost of existing cables and apparatus is included, to analyse the influence of the historical performance of the electric items on the investment decisions. A multihorizon methodology is developed to take into account the long term variation of the demand, combined with the long term variation of cables maintenance costs. Moreover, technological learning coefficients are considered, to take into account the investment costs reductions that arise when an investment in network restructuring and/or reconfiguration is repeated throughout the years. Finally, construction time constraints are included to find a proper investment scheduling that allows a feasible power flow, also during the years required to build a new connection or restructure an existing one. 
This study is also providing recommendations for future research directions within the power system reliability field. The analyses show the important and urgent need of proper methodologies for a better definition of cables projected maintenance costs and learning coefficients dedicated to network restructuring, reconfiguration and expansion.
\end{abstract}

% keywords can be removed

\section{Introduction}\label{intro}

{\renewcommand{\arraystretch}{1.4}
\begin{table}[!htbp]
\label{tab:indexes}
%\scriptsize

\begin{center}
\begin{tabular}{ll}
\bfseries Nomenclature & \\
\hline
\textbf{Indexes}\\
\hline
t & Time step\\
a & Years \\
i,j & Nodes of the grid\\
w & Renewable power plants\\
g & Conventional power generators\\
\hline
\bfseries Parameters & \\
\hline
$D_{i,t,a}$ & Power demand\\
$r$ & Interest rate\\
$G^{op}_g$ & Operational cost of the conventional generator\\
%G^{min}_g & Minimum power out of the conventional generator\\
$G^{cap}_{g,t}$ & Capacity of the conventional generator\\
$W^{cap}_{w,i}$ & Capacity of the renewable plant\\
$W^{\%}_{i,t}$ & Percentage of renewable generation that can be produced in each node \\
%$W^{max}_{w,i,t}$ & Maximum production of each renewable plant\\
$E_{i,j}$ & Binary parameters equal to 1 if a line exists between node i j\\
$\overline E_{i,j}$ & Capacity of the cable\\
$Z$ & Construction time \\
$BigM$ & A very big number\\
$N^{pot}_{i,j}$ & Parameter that is equal to 1 if a potential arc can be placed between node i and j\\
$N^{cost}_{c}$ & Cost of a new cable of type c\\
$N^{life}_{c}$ & Life of a new cable of type c\\
$N^{life}_{i,j}$ & Residual life of an existing cable between nodes i and j\\
$\overline N_c$ & Capacity of the new cable of type c\\
$E^{maint}_{i,j,a}$ & Maintenance cost of existing cables\\
$X_{i,j}$ & Binary parameter equal to 1 if the users wants to evaluate replacement of an existing cable, 0 otherwise\\
$L$ & Technological learning coefficient\\
$\beta_{i,j}$ & Reactance of cable between nodes i and j\\
\hline
\bfseries Variables & \\
\hline
$f^{conv}_{g,i,t,a}$ & Power produced from the generator\\
$f^{ren}_{w,i,t,a}$ & Power produced from the renewable plant\\
$p_{i,j,t,a}$ & Power flow in each arc \\
$d_{i,j,t,a}$ & Binary variable equal to 1 if the flow is into the node i and 0 otherwise\\
$y_{i,j,c,a}$ & Binary variable equal to 1 if a potential arc is created between nodes i and j, 0 otherwise\\
$k_{i,j,c,a}$& Binary variable equal to 1 if an existing arc i,j is replaced by a potential cable of type c, 0 otherwise\\
$m_{i,j,a}$ & Calculated maintenance cost of an existing cable in each year\\
$k^P_{a}$ & Technological learning for new potential cables installation activities\\
$k^R_{a}$ & Technological learning for restructuring activities\\
$q^P_{i,j,c,a}$ & Cost reduction for new cable installation activities due to technological learning\\
$q^R_{i,j,c,a}$ & Cost reduction for restructuring activities due to technological learning\\
$\theta_i$ & Voltage angle on node i\\
\hline
\end{tabular}
\end{center}
\end{table}
}

The electrical power demand is increasing in twice the rate of the overall energy demand \cite{WEO-IEA}. The European Union renewable energy directives assigned the target of at least 27\% renewable energy share in the total energy consumption by 2030 \cite{EESEU}. The share of non-dispatchable and intermittent resources such as wind and solar power has also significantly increased in the last decade. Consequently, the electrical power system is facing challenges due to the requirement of ensuring the quantity and quality of power on demand instantaneously. In fact, network congestion arise due to existing capacity limitations, and the up-front investments in the network make it harder to find a feasible solution. As consequence, optimal capacity expansion of the transmission lines and total generation capacity is needed. In addition, local electricity markets together with novel peer to peer communication and trading mechanisms between microgrids, are changing the traditional power system paradigm \cite{cornelusse2019community}, \cite{morstyn2018using}. Peer to peer electricity trading is still at early stages for microgrids and local communities, however it can affect the way investments in generation expansion are taken and it can add new challenges in the field of reliability and adequacy of the power system.

%Moreover, the electric utility infrastructure is undergoing digitization. Both production and asset information are now widely available for decision makers. Hence, investigating modern power system planning is even more important.
The novel peer to peer electricity trading between microgrids, can result in new physical connections that need to be built among the microgrids, in place of (or in addition to) new investments in conventional or renewable sources \cite{mishra2019multi}.
However, expansion decisions made for an area of the network may affect the power quality and power system reliability in other areas of the power network \cite{mishra2019rnr}. The power quality is the ability to supply power by keeping the electrical parameters (i.e. frequency, voltage, current) within the pre-set limits \cite{wang2010modern}. The power system reliability refers to the ability of a power system to meet customers' requirements for electrical energy and to provide adequate, secure and stable power flow in a given network. The general area of reliability is usually divided into the key aspects of system adequacy, system security and system stability. The system adequacy refers to the existence of sufficient facilities within the system to satisfy customers demand (namely, facilities to generate sufficient energy as well as transmission and distribution facilities to transport energy to the load) \cite{billinton1994basic}. The system security refers to the ability of the system to respond to disturbances arising in the system (i.e. generation/line outages) \cite{morison2004power}. The system stability (frequency, voltage, transient) is the ability of the system to return to its normal or stable conditions after being disturbed.
Power network reliability indicators can be broadly classified into four categories: (a) life-cycle of equipment, (b) environmental factors, (c) consumers experience, (d) fault clearance, and system maintenance \cite{vcepin2011assessment}, \cite{mishra2019reliability}.

Therefore, the profit maximisation objectives that derive from the capacity expansion, must be reached by fulfilling the technical constraints that derive from the necessity of keeping the overall power system quality and reliability indicators within proper boundaries, in order to maintain the network in feasible conditions. 
Keeping the overall feasibility of the network within proper boundaries, will require additional investments in network reinforcement (i.e. increase capacity of existing transmission lines) and reconfiguration (i.e. remove existing cables and/or build new cables somewhere else) also in those areas that reside much farther from the actual nodes where the capacity expansion is actually happening. Such additional investments have to be considered as additional costs when maximising the capacity expansion profits. 
A holistic techno-economic approach is therefore needed and crucial to understand how investments in expansion, reinforcement and reconfiguration change when all the main dimensions of the power system quality and reliability are considered. 

\subsection{The Challenge of Cables Maintenance Costs}
When planning investments in network reinforcement and reconfiguration, two main cost factors must be considered: not only the cost of purchasing and installing new cables and apparatus, but also the future costs of keeping the existing cables, substations, and electrical apparatus in the current operational state. The former is straightforward to calculate, as it relates to the capital cost of cables and apparatus available in the market, together with installation/construction costs \cite{conejo2018power}. The latter is tricky, as it refers to a representative maintenance cost of existing cables and apparatus. This cost differs for each single case, according to the historical performance of the specific electrical item, which depends on many factors (i.e. number of faults experienced, age, number of maintenance interventions etc). For example, a cable that is old and has a history of many faults, is likely to require higher maintenance costs in the forthcoming future, compared to a cable that historically did not experience any critical fault. Therefore it may be economically beneficial to invest in the first cable replacement, and keep the second one as it is for some more years, while continuously monitoring its performance. The investment costs for purchasing and installing a new cable, might be lower than the future costs of keeping it, fixing its failures, and paying for maintenance interventions.

\subsection{The Challenge of Network Feasibility and Construction Time}
When planning investments in network reinforcement and reconfiguration, the construction time will play an important role \cite{zhang2018candidate}. When one or more cables must be reinforced or replaced, the power flow in that corridor is obstructed, but the overall system will have to satisfy the load by finding alternative and temporary solutions. Therefore a multi-horizon perspective to plan investments over many years ahead and guarantee network feasibility throughout the construction time, becomes very important. Cables reinforcement and replacement cannot always happen simultaneously: they have to be planned sequentially in such a way that the overall power system quality and reliability will always be guaranteed.

\subsection{The Novel Role of Technological Learning}
Technological learning effects can play an important role, when defining investments not only here and now, but also scheduling them over a longer period of time of several years. Learning curves express the fact that experience is required, if a technological choice is going to improve and become competitive. That is, technologies will not evolve unless experience with them is possible \cite{barreto2001technological}. This concept can have important implications when it comes to cables replacements and reinforcements made along different corridors of a network. Due to learning effects, an investment will have a lower cost due to the experience gained if the same investment have been already performed in previous years. When evaluating upgrades within different corridors of a network, the combinations of learning effects together with construction time requirements, as well as future demand projections, can affect the scheduling of the investments.

\subsection{Paper Objectives}
The main objective of this paper is to present a techno-economic optimisation tool in order to: 1) study how the power system expansion decisions affect the reliability and adequacy of the network; 2) define how the power system expansion decisions can be taken in a more economical and efficient way by minimising the consequent costs of network reinforcement and reconfiguration; 3) investigate how the network reinforcement and reconfiguration actions should be planned within a long term time horizon of several years, by continuously keeping the network feasibility and ability to satisfy the load.

\section{Key Contribution}
The key contribution of the proposed paper is both methodological and analytical. On the methodological side, this paper represents an extension of the Reliability oriented Network Restructuring (RNR) problem discussed in \cite{mishra2019rnr}. New modelling approaches for the reliability-oriented network restructuring problem are developed and tested, by including key additional features that are currently not addressed in literature, and that enhance the decision making process. 
In particular, the methodological key contributions of this paper can be summarised as follows:

\begin{itemize}

\item Inclusion of maintenance cost of existing cables and apparatus in order to analyse the influence of the historical performance of the electric items on the investment decisions;

\item Multihorizon methodology to take into account the long term variation of the demand combined with the long term variation of cables maintenance costs, for the investment decision making in grid reconfiguration, restructuring and expansion.

\item Inclusion of construction time in the optimisation model, in order to find a proper investment scheduling that allows a feasible power flow also during the years required to build a new connection or restructure an existing one.

\item Inclusion of technological learning coefficients in the optimisation model, to take into account the investment costs reductions that arise when an investment in network restructuring and/or reconfiguration is repeated throughout the years.
\end{itemize}

On the analytical side, computational experiments are performed to validate the tool, and show the impact of the different features on the investment decision making process. In particular, the following analyses are proposed:

\begin{itemize}

\item Extensive sensitivity analyses to show the consequences of technological learning effects on the investment decision making and on the investments scheduling. 

\item Extensive sensitivity analyses to show the combined consequences of future projections of demand trends, together with future projections of cables maintenance costs. 

\item Extensive sensitivity analyses to discuss the trade off between new cables installation, cables restructuring, and cables dismantling in light of the future projections of demand, cables maintenance costs and learning effects.
\end{itemize}

The rest of this paper is organised as follows. Section \ref{LitReview} will provide an overview of the scientific literature in the field of power systems reliability as well as an overview of the available literature in the field of multihorizon modelling for power systems design. The following Section \ref{Model} will discuss the proposed mathematical model with a particular focus on the mathematical modelling of novel features for the reliability oriented network restructuring problem. Computational experiments will be discussed in Section \ref{Analyses}, while conclusions and recommendations for future research directions will be drawn in the last Section \ref{Conclusions}.

%%%%%%%%%%%%%%%%%%%%%%%%%%%%%%%%%%%%%%%%%%%%%%%%%%%%%%%%%%%%%%%%%%%%%%%%
\section{Literature Review}\label{LitReview}

Even though quality and reliability are one of the major concerns in power system and have received wide attention in literature, a holistic techno-economic perspective like the one outlined in the previous section, is lacking in most of the available studies. A holistic techno-economic perspective should take into account the system requirements as a whole: it should handle the network expansion investment decisions by considering also the reinforcement and reconfiguration requirements, together with the maintenance cost of existing electrical apparatus, construction time constraints, and technological learning opportunities. This should be done using a multi-horizon approach for investment planning and a proper methodology to schedule investments in new electrical apparatus throughout the years. The investment scheduling should ensure network feasibility and network security at all times. It should also take into account the possibility to learn from experience in order to make future updates of the electrical apparatus more economical.
%prioritise the most important investments within the most critical areas.
Utilities have two main objectives: 1) get maximum benefit from power apparatus by better utilising the assets they own in a way that is prolonging the remaining utilisation life; 2) maintain the operating condition of the power network. Most of the models and tools available in literature, are developed with the objective of determining optimal asset management strategy while effectively utilizing the collected network information. They usually focus on just one of the aspects mentioned in the previous paragraphs, by overlooking a more holistic techno-economic perspective. In [13] a tool for distribution network reliability analysis is presented. In [14] the factors affecting failure rate in the power system are outlined. In [15] a probabilistic model for evaluating the reliability in a distribution network is presented. In [16] an optimization model is presented that takes into account outage costs and costs of switching devices, along with the nonlinear costs of investment, maintenance and energy losses of both the substations and the feeders. In [17] reliability analysis of a composite power system with renewable resource that is wind farm is presented. In [18] the issues related to reliability, economic and environment for a microgrid with high share of renewable are studied. A big data oriented asset management for electric utilities is presented in [19].
%Traditional failure rate models in use for reliability analysis are based on a constant failure rate without any consideration of operational stress, planning horizon, social implications, environmental factors of the network and components [20-21]. In addition to that, solving a complex model that considers both asset and network information is less efficient and computationally demanding. Many of the indicators such as consumer satisfaction index, fault duration and fault severity change over time, thus more sophisticated machine learning approaches coupled with techno-economic mathematical optimisation models are necessary to address the different aspects of the problem in a holistic way. 

%%%%%%%%%%%%%%%%%%%%%%%%%%%%%%%%%%%%%%%%%%%%%%%%%%%%%%%%%%%%%%
In \cite{baran1989network, shirmohammadi1989reconfiguration, peponis1995distribution,rao2012power,choi2000network,kashem1999network,syahputra2014optimal}, the impact of network reconfiguration in distribution system in terms of loss reduction and load balancing is presented. Authors in \cite{zhu2002optimal,su2003network} present an evolutionary genetic algorithm and mixed-integer hybrid differential evolution based distribution network reconfiguration. A path based reconfiguration technique is presented in \cite{ramos2005path}. How service restoration can be achieved by network reconfiguration is discussed in \cite{shirmohammadi1992service}. Assessing the importance of node in a power network is discussed in \cite{liu2007node}. In \cite{amanulla2012reconfiguration} the power network reconfiguration issue was addressed considering reliability and power losses. In \cite{lopez2004online} the reconfiguration based on variable demand is presented. The study in \cite{zhang2012reliability} presents a network reconfiguration considering the uncertainties of data by interval analysis. In \cite{jose2016reliability,roytelman1995multi} a distribution feeder reconfiguration considering the reliability to minimize the power loss is presented. Adequacy assessment of generating systems with wind power was studied in \cite{liu2017network,gao2009adequacy, billinton2004generating}. Time-series models for reliability evaluation of power systems containing wind energy systems are presented in \cite{billinton1996time}. In \cite{xu2012adequacy} the authors studied the feasibility of distributed energy systems containing electric energy storage and renewable resources adequacy and economic aspects. A probability based reliability assessment is presented in \cite{maruejouls2004practical,gu2012reliability}. A robust power distribution network planning model is presented in \cite{dehghan2015reliability}. A two stage power network reconfiguration strategy considering node importance and restored generation capacity is presented in \cite{zhang2014two}. An evolutionary heuristic technique is applied for skeleton network reconfiguration that is topological characteristics of scale free network in \cite{liu2007skeleton}. A network reconfiguration approach combining genetic algorithm and fuzzy logic is presented in \cite{huang2002enhanced}. A multilevel graph approach for power network reconfiguration is presented in \cite{li2009power}. Energy regulator supply restoration time is studied in \cite{ridzuan2019energy}. An optimal scheduling model for re configurable microgrids is presented in \cite{hemmati2018optimal}.

%%%%%%%%%%%%%%%%%%%%%%%%%%%%%%%%%%%%%%%%%%%%%%%%%%%%%%%%%%%%%%%

As discussed above, there have been many studies conducted on the reconfiguration of the power network. However the distinction between where the re-organization takes place and how does it change the objective is not investigated. For instance there is a significant difference between re-organizing the core network and extending the existing network. The authors of this paper classify the former as a reconfiguration and the latter as a restructuring. %Even though network reliability criterion have been explored in the literature, they are focused on the network as a whole. While the proposed work, aims at associating both the reliability and adequacy attributes to the nodes and arcs of the power network. There is a difference between investigating the reliability and adequacy of the whole network, compared to investigating it for each corridor. Indeed the opportunity to focus on each corridor allows the system operator to identify the reliability issues in more detail. 
The model proposed in this paper, is making the decision of restructuring and reconfiguration holistically by considering at the same time the reliability and adequacy issues.
%Furthermore, the proposed model is considering the ability to improve over time through experience. - learning effect is introduced along with multi-horizon effect for real-world applications. 

Cables maintenance cost definition and future projections are not addressed in literature. The proposed paper aims at showing the role that such feature can play within the investment decision making process for network restructuring and reconfiguration. The performed sensitivity analyses, show that the choice between keeping an existing connection as it is or changing it, is highly dependent on the forecast behaviour of the existing apparatus. 

%The proposed paper aims at discussing the inclusion of a representative maintenance cost within mathematical optimization models for the reliability oriented network restructuring problem. Sensitivity analyses show how this affects the decision-making process.

Technological learning formulations have been widely used in literature and within energy optimisation models, by mainly focusing on the generation mix for CO2 emission reductions and low carbon pathways. However, there are no previous works applying technological learning formulations to the reliability oriented network restructuring problem and to the particular case of cables and electrical apparatus upgrades. In \cite{xu2016bottom} learning coefficients for different low carbon technologies are introduced, and the possibility of equipment upgrading or early decommission of outdated technologies is investigated. The learning rate is introduced to describe the decrease rate of the investment cost when the cumulative experience doubles. Learning rates of different technologies are also included in \cite{rentizelas2012investment}. The model is used to assist long-term investment planning in the electricity production sector and to define the future electricity generation mix, up to the year 2050.
A power capacity expansion problem with technology cost learning is discussed in \cite{heuberger2017power} where two cost learning curves for the different power technologies are derived. The study shows how the inclusion of technological learning is affecting the decisions, the timing of investments, and the competitiveness of technologies. A broad overview of technological learning in energy optimisation models is also provided in \cite{barreto2001technological}, where authors discuss the meaning and impact of learning curves when analysing the transition to low carbon technologies in the energy sector.

All the works mentioned above, apply technological learning coefficients to different renewable and conventional technologies, such as hydro, solar, wind, coal, nuclear, carbon capture and storage etc. To our knowledge this is the first time that technological learning is applied to a reliability network restructuring problem, by considering the learning opportunities gained when upgrading existing electrical apparatus, and scheduling the upgrade of different corridors in the power network.

The multihorizon methodology has been adopted in literature to study some energy related problems. Multihorizon is a modelling framework developed to address long term investment decision making, where investment decisions can be taken in different strategic moments, for example, each year. It allows to take into account the long term variation of various parameters (i.e. future projections of demand, electricity prices, investment costs of resources) within the decision making process.
An introduction to the main concept of multihorizon programming with stochastic implementation can be found in \cite{kaut2014multi}. The main real world applications of multihorizon models found so far in literature are related to natural gas infrastructures \cite{hellemo2013multi}, hydro plants management \cite{abgottspon2016multi}, load management in buildings \cite{rocha2016energy}, and European power systems models for the transition to a low carbon future \cite{skar2014future}.
To our knowledge, none of the above studies considers the problem of reliability and network restructuring with a multihorizon perspective, that includes the trade-off between the construction time, the increasing demand and the forecast increment of cables maintenance costs in the future years.

Moreover, the multihorizon structure of the proposed model allows the inclusion of construction time constraints, which means a better definition of the long-term scheduling of investments in network restructuring and reconfiguration. This approach stands out compared to existing approaches, where the investment decisions are mostly taken here and now, without focusing on the practical long term issues that arise when the actual investment implementation happens, and the overall network feasibility has to be granted in each operational time step.

%%%%%%%%%%%%%%%%%%%%%%%%%%%%%%%%%%%%%%%%%%%%%%%%%%%%%%%%%%%%%%%%%%%%%%%%%%%5
\section{Mathematical Model}\label{Model}
This section will outline the mathematical model developed for the reliability oriented network restructuring problem, with a multihorizon perspective, including technological learning, construction time, and cables maintenance costs. The proposed model is based on a DC power flow, however, OPF modelling is also possible, as thoroughly described in \cite{mishra2019rnr}. Of course OPF presents the advantage of a better representation of the technical properties of the grid, together with the disadvantage of dealing with non linear modelling and therefore having to consider sub-optimal solutions. 
For the purposes of this paper, a DC formulation has been chosen to discuss the model, because the main focus of the sensitivity analyses is to show the implications of technological learning, construction time and cables maintenance costs on the decisions making process, rather than focusing on the implications of particular technical properties of the grid. The same model can anyway be tested with the OPF formulation proposed in \cite{mishra2019rnr} to further ensure the feasibility of the investment decisions. 

%%%%%%%%%%%%%%%%%%%%%%%%%%%%%%%%%%%%%%%%%%%%%%%%%%%%%%%%%%%%%%%%%%%%%%5
\subsection{Objective Function}
\begin{scriptsize}

\begin{equation} \label{eq:OF}
\centering
{\min C^{op} + C^{inv}}
\end{equation}

\begin{equation} \label{eq:oper}
\centering
{C^{op} = \sum_{t,i,g,a} f^{conv}_{g,i,t,a}*G^{op}_g
}
\end{equation}

\begin{equation} \label{eq:inv}
\centering
{C^{inv} = \sum_{i,j,c,a} F^{cab}_c *y_{i,j,c,a}*N^{cost}_{c} - q^P_{i,j,c,a} + \sum_{i,j,c,a} F^{cab}_c * k_{i,j,c,a}*N^{cost}_{c} - q^R_{i,j,c,a} + \sum_{i,j,a} F_{i,j} * m_{i,j,a}}
\end{equation}

\begin{equation} \label{eq:CRF1}
\centering
{F^{cab}_c = \frac {r * (1+r)^{N^{life}_c}} {(1+r)^{N^{life}_c} -1}}
\end{equation}

\begin{equation} \label{eq:CRF2}
\centering
{F_{i,j} = \frac {r * (1+r)^{N^{life}_{i,j}}} {(1+r)^{N^{life}_{i,j}} -1}}
\end{equation}

\end{scriptsize}

The objective function \ref{eq:OF} minimizes the total cost that is comprised of operational plus investment costs. Operational costs in \ref{eq:oper} are related to fuel consumption in dispatchable generation. The investment cost in \ref{eq:inv} is the summation of three terms: the cost of installing new potential cables where a connection does not exist, the cost of replacing existing obsolete cables with new ones, and the cost of keeping existing cables as they are (a so called "maintenance cost"). The latter is a representative cost that incorporates all the costs that a company should face to keep a cable as it is. This cost is defined according to the history of the cables, its maintenance requirements, failures and issues. Section \ref{model-restruct} will discuss the equations needed to include the maintenance cost $m_{i,j,a}$ depending on the investment choices.
All the investment costs are actualised using a capital recovery factor, that is calculated as a function of the interest rate $r$ and the lifetime of cables. For new cables installation, the forecast lifetime of a cable $N^{life}_c$ is used in the calculation of the capital recovery factor $F^{cab}_c$. The last term referring to the maintenance cost, is actualized through the capital recovery factor $F_{i,j}$ by taking into account the forecast residual lifetime $N^{life}_{i,j}$ of the existing cables. 

The first two investment costs are discounted using the variables $q^P_{i,j,c,a}$ and $q^R_{i,j,c,a}$ in order to take into account the technological learning that arise if the same type of activity (new cable installation or restructuring) has been performed in the previous years. Section \ref{model-learn} will discuss the equations needed to define the values of the $q^P_{i,j,c,a}$ and $q^R_{i,j,c,a}$ variables.

%\texttt{ADD THE ZONAL COSTS IN THE OF AND THE DESCRIPTION IN THE TEXT AS WE HAVE IN THE RNR DRAFT?}

%\texttt{ADD EQUATION ABOUT crf}

%%%%%%%%%%%%%%%%%%%%%%%%%%%%%%%%%%%%%%%%%%%%%%%%%%%%%%%%%%%%%%%%%%%%%%5
\subsection{Conventional and Renewable Generators}
\begin{scriptsize}
\begin{equation} \label{eq:conv}
\centering
{f^{conv}_{g,i,t,a}<=G^{cap}_{g,t} \qquad \forall (g,i,t,a)
}
\end{equation}

\begin{equation} \label{eq:ren}
\centering
{f^{ren}_{w,i,t,a} <= W^{cap}_{w,i} * W^{\%}_{i,t} \qquad \forall (w,i,t,a)
}
\end{equation}
\end{scriptsize}

This group of constraints limit the capacity of the dispatchable generators as well as the capacity of the renewable plants.

%%%%%%%%%%%%%%%%%%%%%%%%%%%%%%%%%%%%%%%%%%%%%%%%%%%
\subsection{Grid General Management}
\begin{scriptsize}
\begin{equation} \label{eq:grid01}
\centering
{p_{i,j,t,a} <= \overline E_{i,j}  \qquad \forall (i,j,t,a) \quad | \quad E_{i,j}=1 \quad and \quad X_{i,j}=0
}
\end{equation}

\begin{equation} \label{eq:grid02}
\centering
{p_{i,j,t,a} <= BigM * d_{i,j,t,a} \qquad \forall (i,j,t,a)
}
\end{equation}

\begin{equation} \label{eq:grid03}
\centering
{p_{j,i,t,a} <= BigM * (1-d_{i,j,t,a}) \qquad \forall (i,j,t,a)
}
\end{equation}

\begin{equation} \label{eq:grid04}
\centering
{\sum_g f^{conv}_{g,i,t,a} + \sum_w f^{ren}_{w,i,t,a} - \sum_j p_{i,j,t,a}* E_{i,j} + \sum_j p_{j,i,t,a}* E_{i,j} - \sum_j p_{i,j,t,a}* N^{pot}_{i,j} + \sum_j P_{j,i,t,a}* N^{pot}_{i,j} >= D_{i,t,a}
 \qquad \forall (i,t,a)
}
\end{equation}

\begin{equation} \label{eq:grid05}
\centering
{p_{i,j,t,a} = \frac{\theta_{i,t,a} - \theta_{j,t,a}}{\beta_{i,j,t,a}} \qquad \forall (i,j,t,a)
}
\end{equation}

\end{scriptsize}

%\texttt{NOTE: add equation on DC power flow. Cite our RNR paper to explain that OPF can be also utilised with its advantages (more precise power flow calculations) and disadvantages (non linear model)}

Constraint \ref{eq:grid01} imposes that the power flow should be less than or equal to the cable capacity. Constraints \ref{eq:grid02} and \ref{eq:grid03} impose mutually exclusive power flows along each line. Constraint \ref{eq:grid04} defines the flow balance in each node, such that the power flow into each node has to be equal to the power flow out from the node in each time step. Finally, constraint \ref{eq:grid05} imposes the DC power flow main properties to define the voltage angle in each node according to the values of reactance. The voltage angles are restricted to be within feasible limits. More complex OPF formulations are also possible as implemented in \cite{mishra2019rnr} with the advantage of more precise power flow calculations and the disadvantage of the non linear equations that will have to be handled using ad-hoc algorithms.

%%%%%%%%%%%%%%%%%%%%%%%%%%%%%%%%%%%%%%%%%%%%%%%%%%%%%%%%%%%%%%%%%%%%%
\subsection{Restructuring} \label{model-restruct}
\begin{scriptsize}
\begin{equation} \label{eq:res01}
\centering
{p_{i,j,t,a} <= \overline E_{i,j} * \Big( 1 - \sum_{c,a1=1}^{a1=a-Z-1}  k_{i,j,c,a1} \Big) + \sum_{c,a1=1}^{a1=a-Z-1}  k_{i,j,c,a1} * \overline N_c
 \qquad \forall (i,j,t,a) \quad | \quad E_{i,j}=1 \quad and \quad X_{i,j}=1
}
\end{equation}

\begin{equation} \label{eq:res02}
\centering
{\sum_{a1 = a}^{a1 = a+Z-1} p_{i,j,t,a1} <= BigM * \Big( 1- \sum_c k_{i,j,c,a1} \Big)  \qquad \forall (i,j,t,a) \quad | \quad E_{i,j}=1 \quad and \quad X_{i,j}=1
}
\end{equation}

\begin{equation} \label{eq:res03}
\centering
{\sum_{c,a} k_{i,j,c,a} <= 1 \qquad \forall (i,j)
}
\end{equation}

\begin{equation} \label{eq:res04}
\centering
{m_{i,j,a} = \Big( 1 - \sum_{c,a1=1}^{a1 = a} k_{i,j,c,a1} \Big) * E^{maint}_{i,j,a}
 \qquad \forall (i,j,a)
}
\end{equation}

\end{scriptsize}

\begin{figure}[htbp] 
\begin{center}
\includegraphics[scale=0.8]{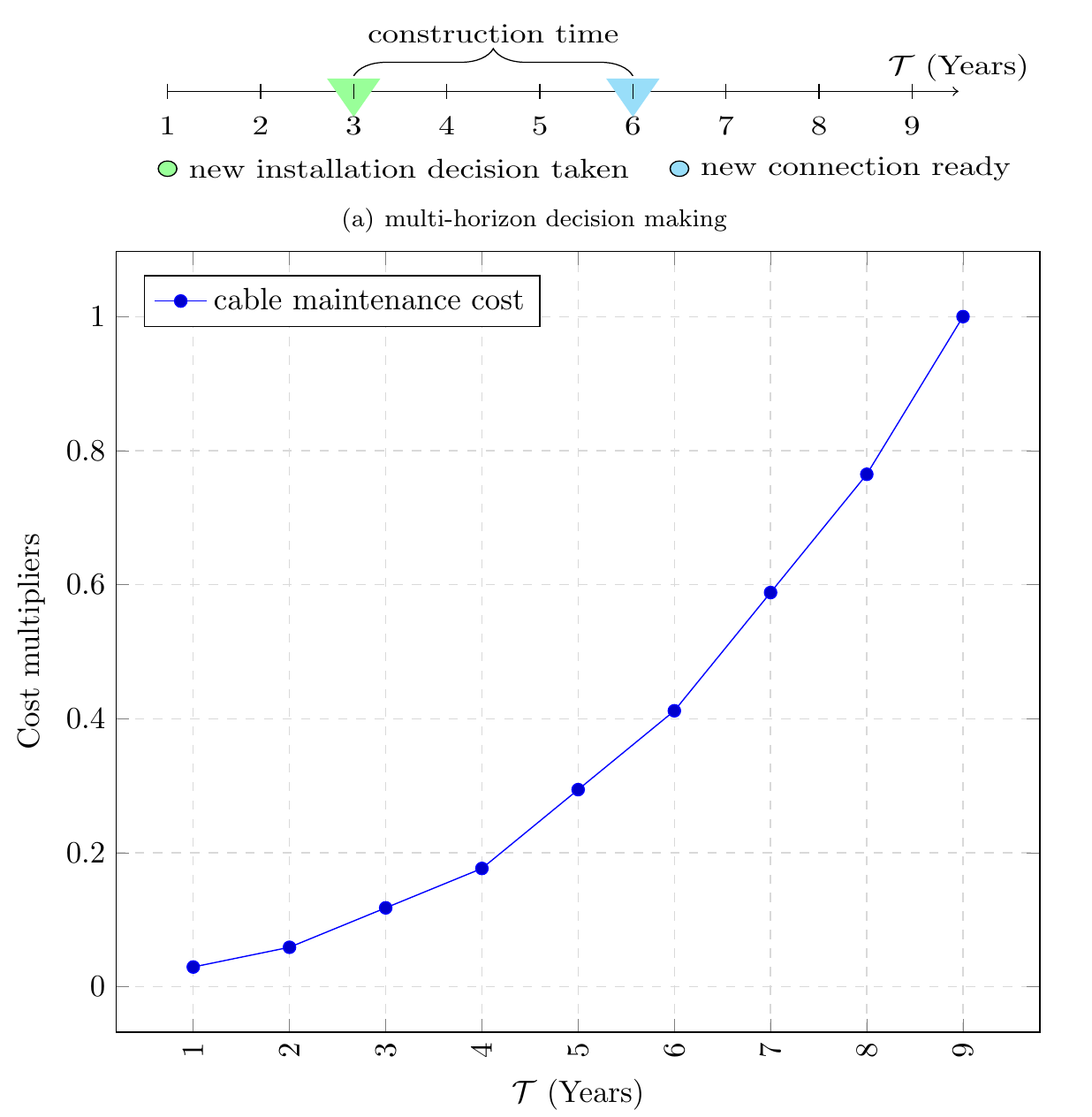}
\caption{Multihorizon structure of the model. When a new installation decision is taken, the power flow will not be allowed throughout the construction time. }\label{Fig-multihor}
\end{center}
\end{figure}

Constraint \ref{eq:res01} defines how reconfiguration can happen. If a cable is replaced with a new one sometimes over the previous years, then the active power in the following year will have as upper bound the new cable capacity. Otherwise the upper bound will be given by the current cable capacity. Construction time is considered, such that the replacement decisions has to be taken in the previous years before the construction starts.
In particular, if restructuring happens, the power flow along the newly installed cables will be allowed only if two conditions are satisfied: the first condition is that sometimes in one of the previous years a replacement decision has been taken; the second condition is that such a decision has to be taken at least $Z$ years before the current year, where $Z$ is the construction time, representing the years required to actually build the new connection. In fact, if a replacement decision is taken in year $a$, it is necessary to consider that $Z$ years will be needed to actually build the new connection: therefore, throughout these construction years, such a connection will not be available and the model will need to ensure that the overall system will still work, even during this transition period, by using alternative feasible corridors within the network. Figure \ref{Fig-multihor} illustrates the multihorizon structure of the model and the trade-off between making a new installation decision at the right time, keeping the system feasibility throughout the construction time, and avoiding increasing maintenance costs of obsolete existing cables. Hence the model has to take decisions not only based on the increased maintenance costs of obsolete cables, but also based on keeping the system feasibility throughout the construction time. In the second sub-figure, an example of the gradual change in the normalized cable maintenance cost over ten years time horizon is presented. This is just an example for illustrative purposes. Different trends for future maintenance cost projections can be considered and tested within the proposed mathematical model.

%\texttt{this figure is not exactly as it should be. The timeline of the section (a) should be the same of the time line i section (b). The two sections are connected. They refer to the same time line. While in the figure the time line of section (a) is different and shorter than the one of section (b)}

Constraint \ref{eq:res02} imposes that, if an existing cable $i-j$ is replaced with a new cable of type $c$, then the active power has to be zero for the years that are needed to actually build the connection (construction time). Indeed, during the years required to build the new cable, the connection between the nodes is down, due to the construction activities that are taking place.

Constraint \ref{eq:res03} limits the choice of new cables to 1. Only one type of cable should be chosen for every replacement. Furthermore, for every couple of nodes, a replacement can happen only once throughout the time horizon. It is therefore assumed that the new cables will last more than the considered time horizon, which is a reasonable assumption given the long cables lifetime. 

Constraint \ref{eq:res04} imposes that, if an existing cable is not replaced, then the maintenance cost must be equal to the current given maintenance cost. On the other hand, if the existing cable is replaced with a new one, then its maintenance cost becomes zero for the years ahead. The maintenance cost is therefore applied only to existing obsolete cables. New cables are assumed to be in ideal conditions and perform well for the forthcoming years.

Restructuring can happen only in those arcs that the operator is willing to check. Not all the arcs of the grid will be subjected to such decision, therefore the binary parameter $X_{i,j}$ is used to select which arcs to check. In real world situations, the arcs that the investor would check are those who appear obsolete, namely those with a high rate of failures in the past years, or high maintenance costs, or technical properties that appear inadequate for the future demand projections (i.e. too low capacity to accommodate a potential forthcoming new district in the area)

It is straightforward that the model allows also the possibility to simply dismantle existing cables by providing the choice to replace an existing obsolete cable with a different one with a fictitious capacity equal to zero. Such a fictitious cable will have a cost equal to zero, or very low, depending on whether the cable will be just abandoned and no longer utilised, or whether it will be actually removed. Environmental costs of disposal can be included in the overall dismantling cost.

%%%%%%%%%%%%%%%%%%%%%%%%%%%%%%%%%%%%%%%%%%%%%%%%%%%%%
\subsection{Potential Grid}
\begin{scriptsize}
\begin{equation} \label{eq:pot01}
\centering
{p_{i,j,t,a} <= \sum_{c,a1=1}^{a1=a-Z-1} y_{i,j,c,a1} * \overline N_c
\qquad \forall (i,j,t,a) \quad | \quad N^{pot}_{i,j}=1
}
\end{equation}

\begin{equation} \label{eq:pot02}
\centering
{\sum_{c,a} y_{i,j,c,a} <= 1 \qquad \forall (i,j)
}
\end{equation}
\end{scriptsize}

Constraint \ref{eq:pot01}  defines how the installation of new cables where no existing connections are available can happen. If a new cable is installed between nodes $i-j$ (where a cable do not exist yet), then the active power between $i-j$ will have as upper bound the capacity of the new selected cable. Construction time is considered, such that the new cable installation decision has to be taken in the previous years, before the construction starts. Therefore, as explained also in the previous section dedicated to restructuring, no power flow is allowed during construction time.

Constraint \ref{eq:pot02} limits the choice of new cables to 1. Only one type of cable should be chosen for every new installation. Furthermore for every couple of nodes, a new installation can happen only once throughout the time horizon. It is therefore assumed that the new cables will last more than the considered time horizon, which is a reasonable assumption given the long cables lifetime.   

New potential installations can happen only in those arcs that the operator is willing to consider. Not all the arcs of the grid will be subjected to such decision, therefore the binary parameter $N^{pot}_{i,j}$ is used to select which arcs should be considered for new potential installations.

%%%%%%%%%%%%%%%%%%%%%%%%%%%%%%%%%%%%%%%%%%%%%%
\subsection{Learning Effect}\label{model-learn}
\begin{scriptsize}

\begin{equation} \label{eq:learn01}
\centering
{k^P_{a} = \sum_{i,j,c,a1=1}^{a1=a-1} y_{i,j,c,a1} * L  \qquad \forall (a)}
\end{equation}

\begin{equation} \label{eq:learn02}
\centering
{k^R_{a} = \sum_{i,j,c,a1=1}^{a1=a-1} k_{i,j,c,a1} * L  \qquad \forall (a)
}
\end{equation}

\begin{equation} \label{eq:learn03}
\centering
{q^P_{i,j,c,a} >= N^{cost}_{c} * k^P_{a} + (y_{i,j,c,a} - 1) * BigM \qquad \forall (i,j,c,a)
}
\end{equation}

\begin{equation} \label{eq:learn04}
\centering
{q^P_{i,j,c,a} <= N^{cost}_{c} * k^P_{a} + (1 - y_{i,j,c,a}) * BigM \qquad \forall (i,j,c,a)
}
\end{equation}

\begin{equation} \label{eq:learn05}
\centering
{q^R_{i,j,c,a}>= N^{cost}_{c} * k^R_{a} + (k_{i,j,c,a}- 1) * BigM \qquad \forall (i,j,c,a)
}
\end{equation}

\begin{equation} \label{eq:learn06}
\centering
{q^R_{i,j,c,a}<= N^{cost}_{c} * k^R_{a} + (1 - k_{i,j,c,a}) * BigM \qquad \forall (i,j,c,a)
}
\end{equation}

\begin{equation} \label{eq:learn07}
\centering
{q^P_{i,j,c,a} <= y_{i,j,c,a} * BigM \qquad \forall (i,j,c,a)
}
\end{equation}

\begin{equation} \label{eq:learn08}
\centering
{ q^R_{i,j,c,a}<= k_{i,j,c,a}* BigM \qquad \forall (i,j,c,a)
}
\end{equation}

\end{scriptsize}

This set of constraints define how technological learning can be taken into account within the reliability oriented network restructuring problem.

Constraint \ref{eq:learn01} imposes that, if a new cable installation is decided in a certain year, and a cable installation in the previous years happened already, the variable $k^P_{a}$ will be equal to a percentage value $L$. This value will define the learning coefficient and it will be used to make the cable cost for the considered year cheaper, due to the technological learning that happened years before.
The higher the number of installations in the previous years, the higher the percentage value $L$ that will be applied to reduce the current cable cost.
If no new installations are decided in a certain year, then the variable $k^P_{a}$ will be equal to zero.

Constraint \ref{eq:learn02} works the same as the previous ones, but applies to restructuring activities, by utilising the variable $k^R_{a}$.

Constraints \ref{eq:learn03} and \ref{eq:learn04} define the value of the variable $q^P_{i,j,c,a}$ that will be discounted by the cable cost in the objective function. If  $k^P_{a}$ is greater than zero, and a new potential installation happens in a certain year, (variable $y_{i,j,c,a} = 1$), then the new cables investment cost for that particular year will be cheaper due to the technological learning effect. 

Constraints \ref{eq:learn05} and \ref{eq:learn06} work the same as the previous ones, but they refer to the restructuring activities and related learning effects.

Constraint \ref{eq:learn07} imposes that, if a new potential installation does not happen in a certain year, then the related $q^P_{i,j,c,a}$ variable for this kind of activity in that year has to be zero.

Constraint \ref{eq:learn08} imposes that, if a restructuring does not happen in a certain year, then the related $q^R_{i,j,c,a}$ variable for this type of activity in that year has to be zero.

\section{Computational Experiments}\label{Analyses}

Computational experiments have been run in order to investigate how the inclusion of technological learning and maintenance costs of cables affects the decision making processes. Additional analyses have been carried out in order to investigate the trade off between new potential installation decisions, restructuring decisions with cables replacements, and reconfiguration decisions that include cables dismantling.

\begin{figure}[htbp] 
\begin{center}
\includegraphics[scale=0.7]{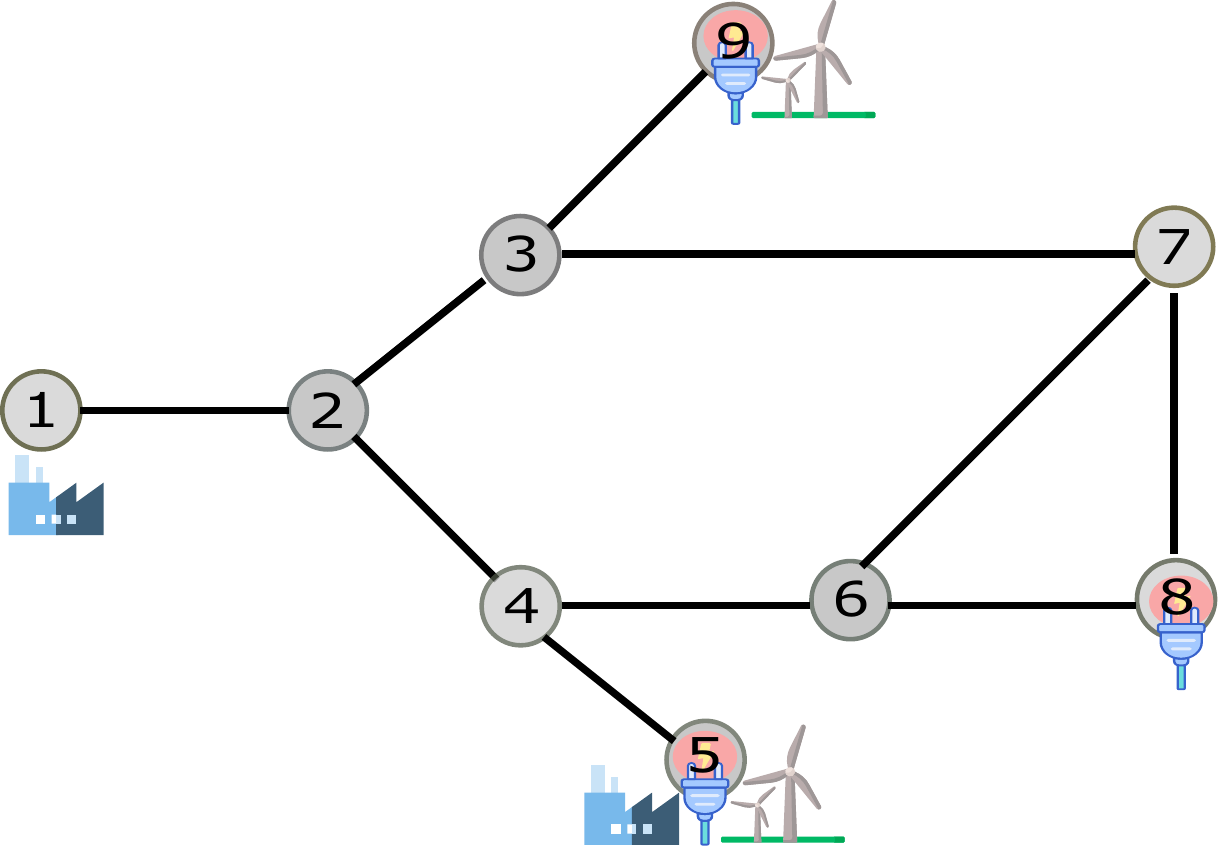}
\caption{Configuration of the network considered for testing}\label{Fig-Slide1}
\end{center}
\end{figure}

The Figure \ref{Fig-Slide1} shows the configuration of the network considered for testing purposes. From a higher level perspective, each arc that is represented in the figure, can also be considered as an aggregated representation of a wider corridor. Therefore the simplified representation of Figure \ref{Fig-Slide1} can also be used to refer to more complex structures, by aggregating corridors according to some zonal properties and clustering techniques.

A time horizon of 10 years has been assumed. Such a period is long enough to observe new cables installation decisions throughout the years, as well as include long term demand variations. An average construction time of three years has been assumed for testing purposes.

\begin{figure}[htbp] 
\begin{center}
\includegraphics[width=0.5\textwidth]{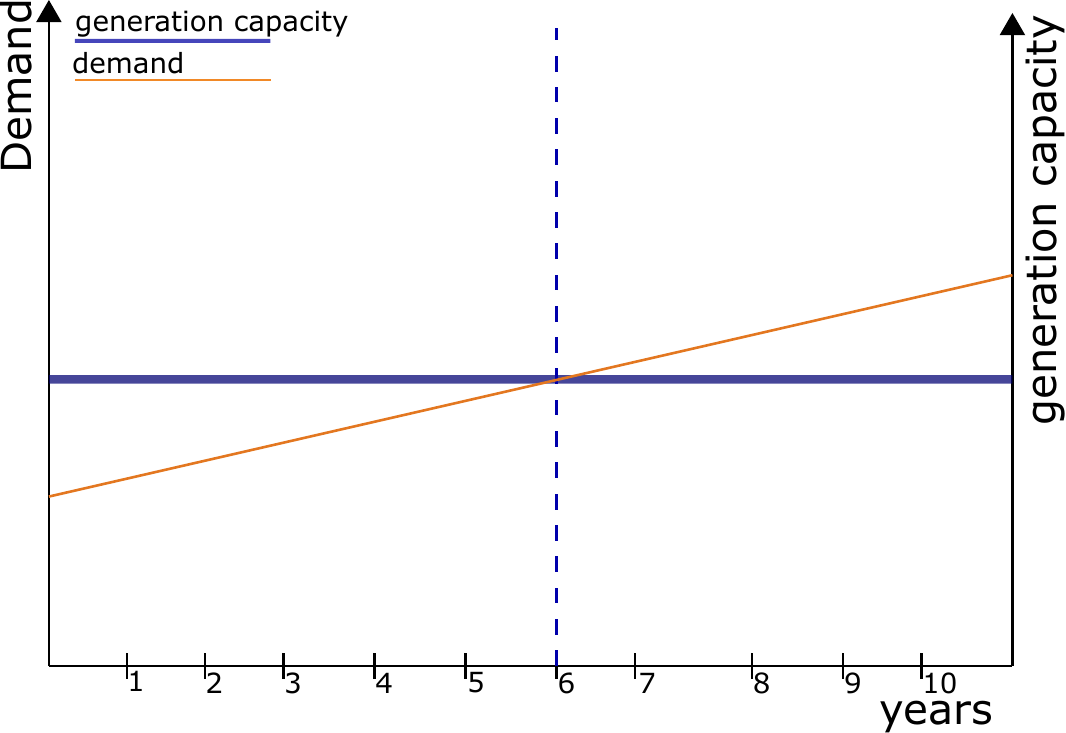}
\caption{Assumptions for node 5. The demand in node 5 will exceed the generation capacity installed in the node from the 6th year}\label{Fig-Node5}
\end{center}
\end{figure}

\begin{figure}[htbp] 
\begin{center}
\includegraphics[width=0.5\textwidth]{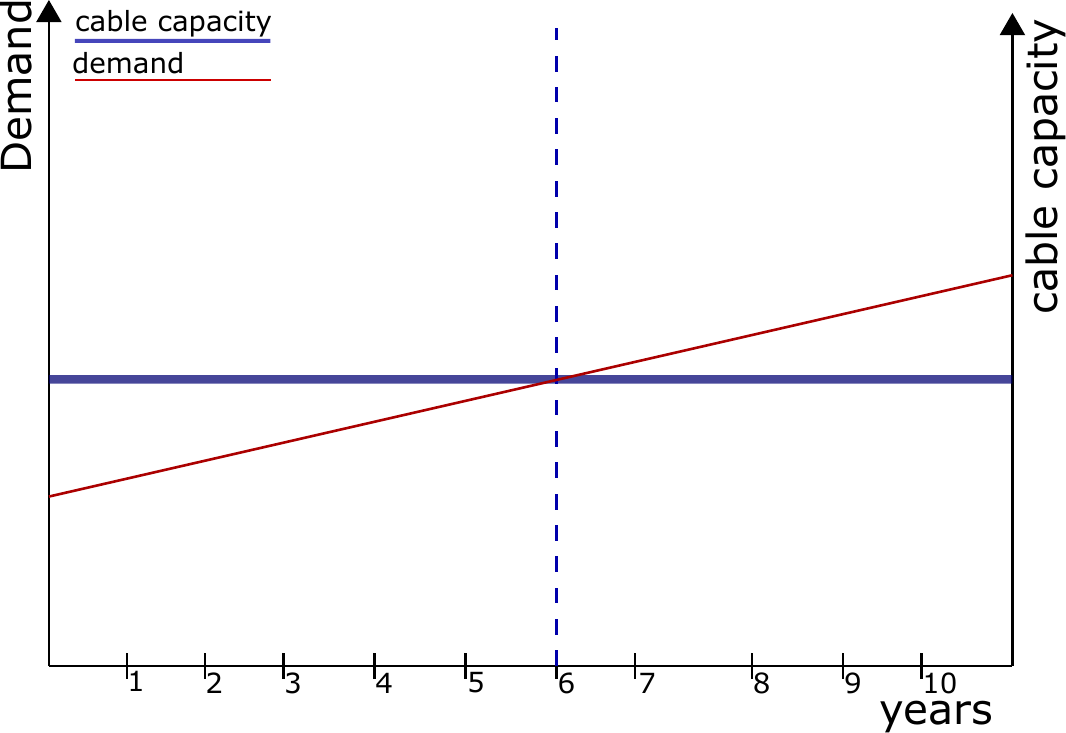}
\caption{Assumptions for node 8. The demand in node 8 will exceed the capacity of the connected cables from the 6th year}\label{Fig-Node8}
\end{center}
\end{figure}

\begin{figure}[htbp] 
\begin{center}
\includegraphics[width=0.5\textwidth]{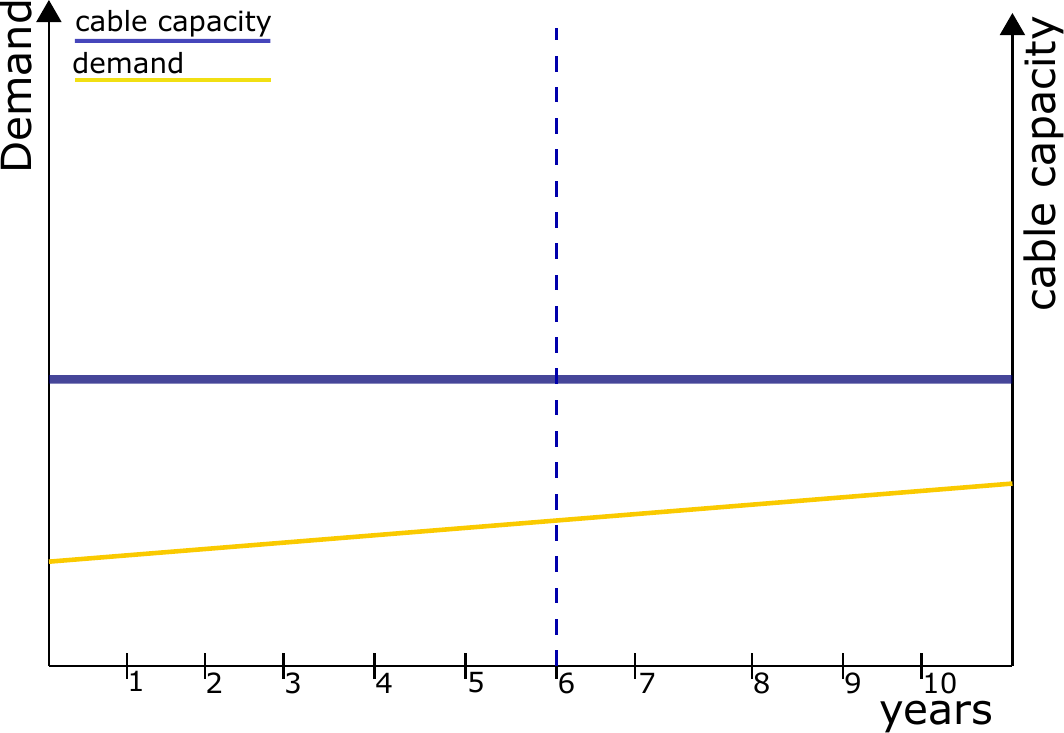}
\caption{Assumptions for node 9. The demand in node 9 will increase very slowly without exceeding any cable capacity}\label{Fig-Node9}
\end{center}
\end{figure}

Node 1 is equipped with a conventional generator whose capacity is assumed to be big enough to satisfy the load of the grid. Node 5 is equipped with a wind plant and a small conventional generator. The demand in node 5 is assumed to increase linearly throughout the years, in such a way that it will exceed the generation capacity installed in the node from the 6$^{th}$ year ongoing. 
Node 8 is assumed to be a district with a demand that is supposed to increase linearly throughout the years, in such a way that it will exceed the capacity of the connected cables from the 6th years ongoing. Finally, node 9 is equipped with a renewable plant and the demand is assumed to be almost constant, with a very low increment throughout the years that will not exceed the available existing generation and cable capacity.
Figures \ref{Fig-Node5}, \ref{Fig-Node8} and \ref{Fig-Node9} summarise the main assumptions described above. 

 The following paragraphs will illustrate sensitivity analyses by focusing the attention on different model features each time. Each figure will show the current configuration of the network (on the left), and how the network configuration changes after the optimisation (on the right). Looking at the left side of each figure, the following symbols and colours will be used: existing cables are represented by black lines; potential new installations between nodes are represented by red dot lines; potential restructuring of existing cables with new ones is represented by a red straight line. When maintenance cost of cables is considered in the experiment, it will be represented by the symbol of a dollar. Such a symbol will have a small size or a big size to identify a cable with either a low or high maintenance cost respectively.
 Looking at the right side of each figure, the following symbols and colours will be used: dismantling decisions are represented by a blue dot line, new installation decisions are represented by a green line, restructuring decisions are represented by a yellow line. An $"Y"$ followed by a number between 1 and 10 will indicate in which year a restructuring, dismantling or new installation decision is taken. For the new installation and restructuring decisions, the corresponding lines will be thin or thick depending on whether the new cable has a big size (and high cost) or a small size (and a low cost) respectively. This will facilitate the figure reading without having to go too deep into numerical details. Indeed, in this kind of analyses, it is not the actual numbers the matters, but rather the ratio between them and the interplay between the different entities. If only one cable size is considered for a restructuring or new installation decision, this will be represented by a single green or yellow line without differentiating between thin or thick.
 
%%%%%%%%%%%%%%%%%%%%%%%%%%%%%%%%%%%%%%%%%%%%%%%%%%%%%%%%%%%%%%%%%%%%%%%%%%%%%
\subsection{Focus - Learning Effect and Cables Maintenance Costs}

\begin{figure}[htbp] 
\begin{center}
\includegraphics[scale=0.7]{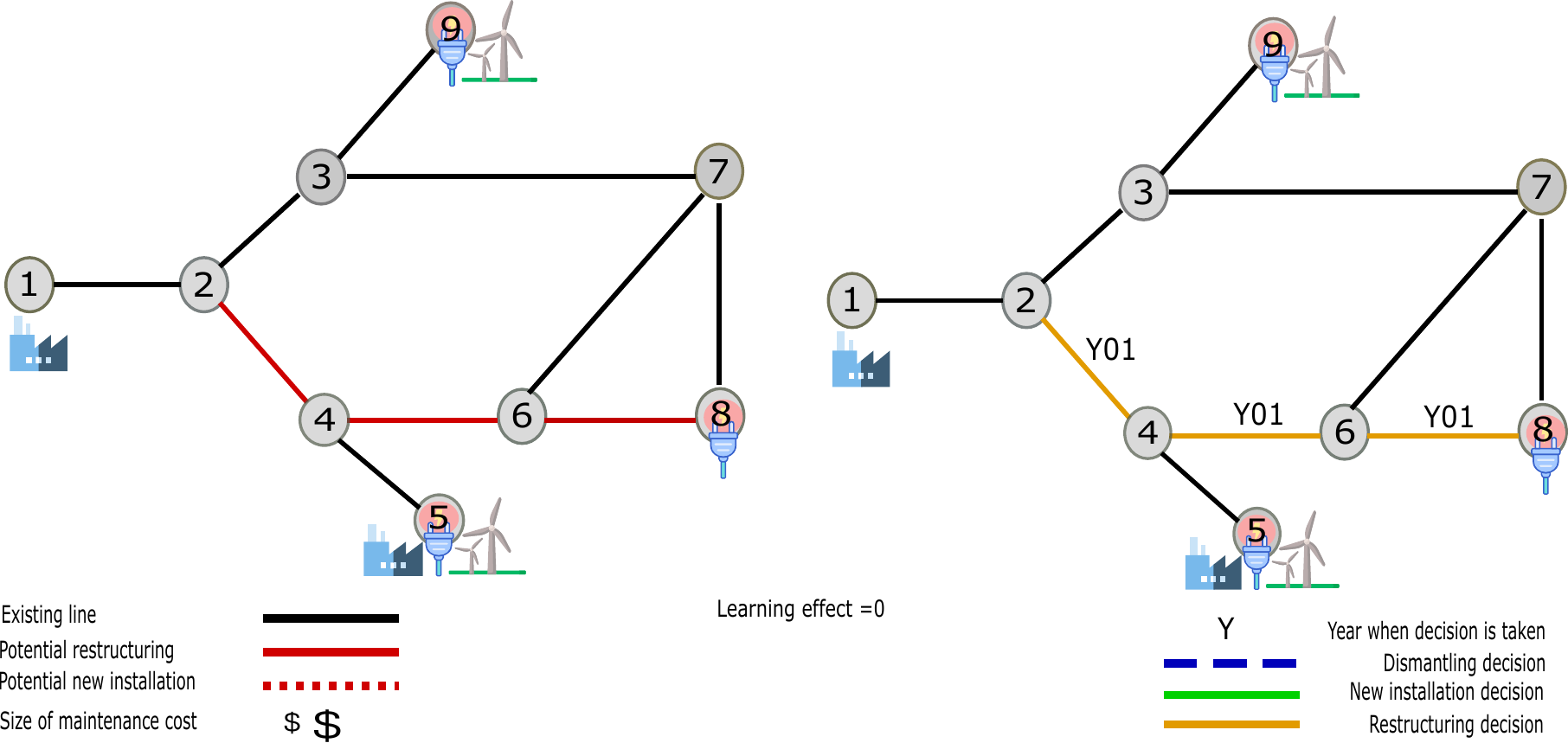}
\caption{Restructuring decision when technological learning is not considered}\label{Fig-Slide2}
\end{center}
\end{figure}

Figure \ref{Fig-Slide2} shows the investment decisions when technological learning is not considered. Cables 2-4, 4-6 and 6-8 are undersized compared to the demand increment that is assumed in the following years, as outlined in the previous paragraph. Therefore a restructuring decision is taken on the first year to replace cables 2-4, 4-6 and 6-8 with bigger ones that will be ready in time to accommodate the higher future power flow. During the construction time, cables 2-4, 4-6 and 6-8 will be down at the same time, but both network feasibility and network security will be guaranteed. Indeed, network feasibility is guaranteed because alternative existing corridors can be used to continuously fulfill the demand in the nodes, while the corridor 2-4, 4-6 and 6-8 will be down. Network security is guaranteed as well, because all the nodes involved by the dismantling tasks, are anyway connected to a generator. In particular, node 5 is directly connected to a conventional generator that is installed in the node itself. As for nodes 8 and 9, they are indirectly connected through arcs to the bigger conventional generator that is installed on node 1. 

%\texttt{QUESTION: should we add a reference about this rule of having at least one generator connected to ensure security?}

%\texttt{QUESTION: I should probably introduce the concepts of network feasibility and security at the beginning}

\begin{figure}[htbp] 
\begin{center}
\includegraphics[scale=0.70]{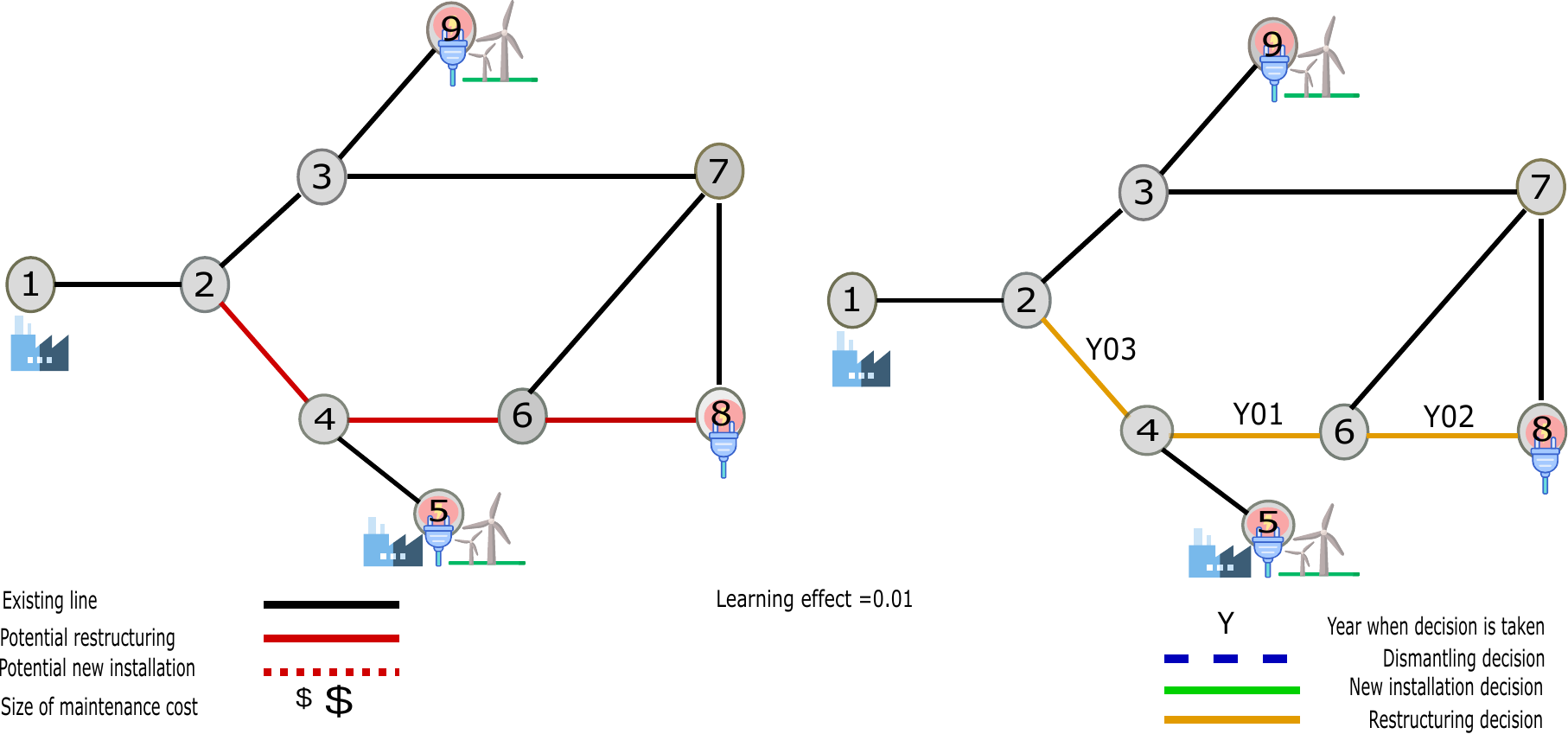}
\caption{Restructuring decision when technological learning is included}\label{Fig-Slide3}
\end{center}
\end{figure}

The Figure \ref{Fig-Slide3} shows how the decision changes when technological learning is taken into account. The concept of technological learning is that, if an investment has been performed sometimes before, repeating the investment again in the future will be cheaper due to the improved knowledge acquired through the past experience. In the particular case of Figure \ref{Fig-Slide3} it has been assumed a learning effect of 1\%, meaning that the restructuring investment will be 1\% cheaper if it has been performed before. The implication of including such a learning effect, is that investments in cables restructuring are now scheduled to take advantage of the experience that can be gained year by year. In figure \ref{Fig-Slide2}, in absence of learning effects, the whole corridor was dismantled at the same time. While in the case of figure \ref{Fig-Slide3}, the possibility to learn from experience, leads to an investment scheduling such that cables are now dismantled in different years. In this case the scheduling takes into account both the learning possibility, the construction time, the need to have the corridor ready within the 6th year, as well as the possibility to use alternative paths in those periods of time where cables are down.

\begin{figure}[htbp] 
\begin{center}
\includegraphics[scale=0.7]{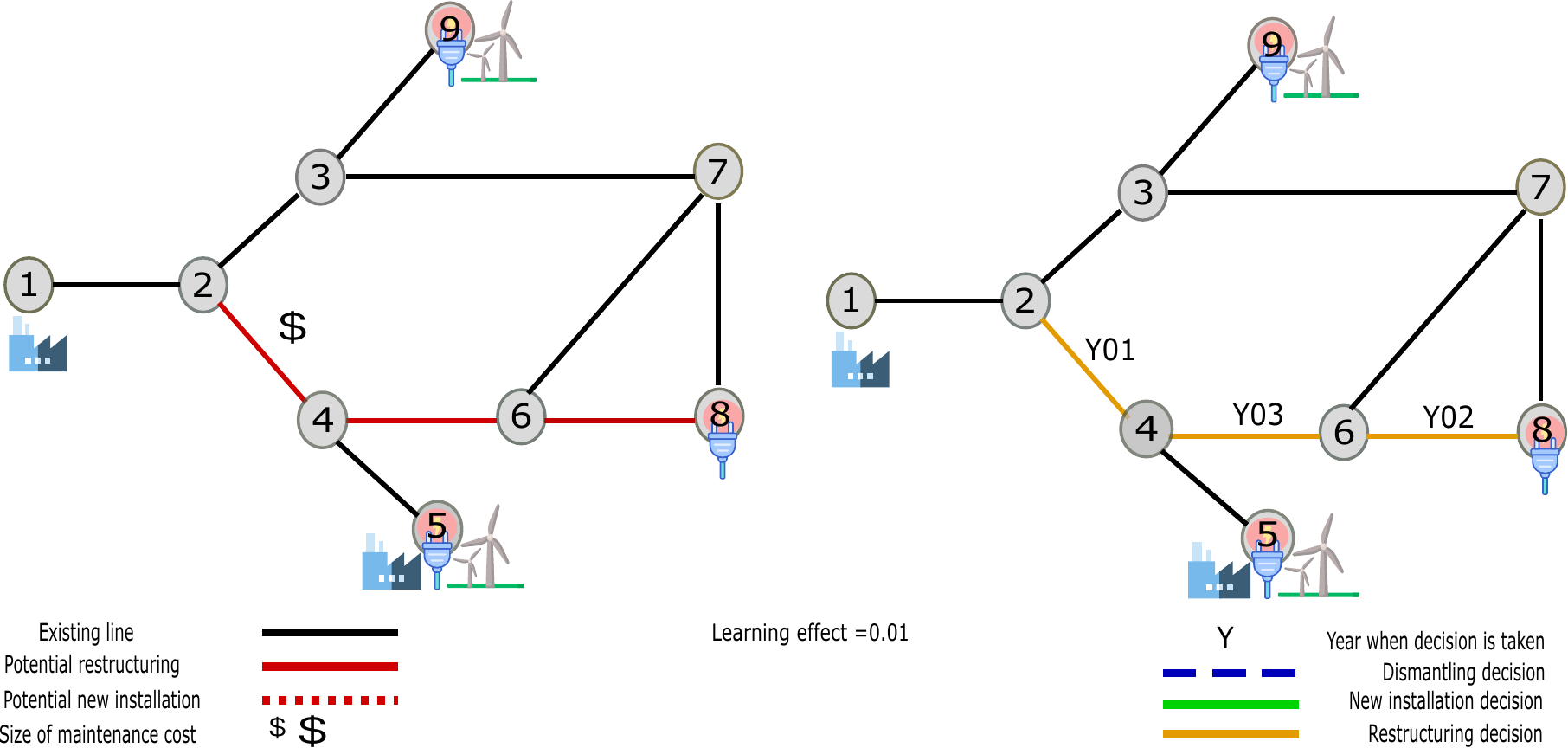}
\caption{Restructuring decision when technological learning and the maintenance cost of one cable are included}\label{Fig-Slide4}
\end{center}
\end{figure}

Figure \ref{Fig-Slide4} shows how the decision changes when also maintenance costs of cables are involved. In this case, it is assumed that cable 2-4 is very obsolete and presents a maintenance cost, compared to the others that are assumed to be still in a good shape. Compared to Figure \ref{Fig-Slide3}, where cable 2-4 was dismantled only on the 3rd year, in Figure \ref{Fig-Slide4} the cable 2-4 is now scheduled for dismantling on the very first year. Indeed, the presence of a maintenance cost is making the cable very expensive and therefore it is better to get rid of it sooner than the others. Therefore, not only learning effect but also maintenance cost affect the scheduling of the investment in restructuring.

\begin{figure}[htbp] 
\begin{center}
\includegraphics[scale=0.7]{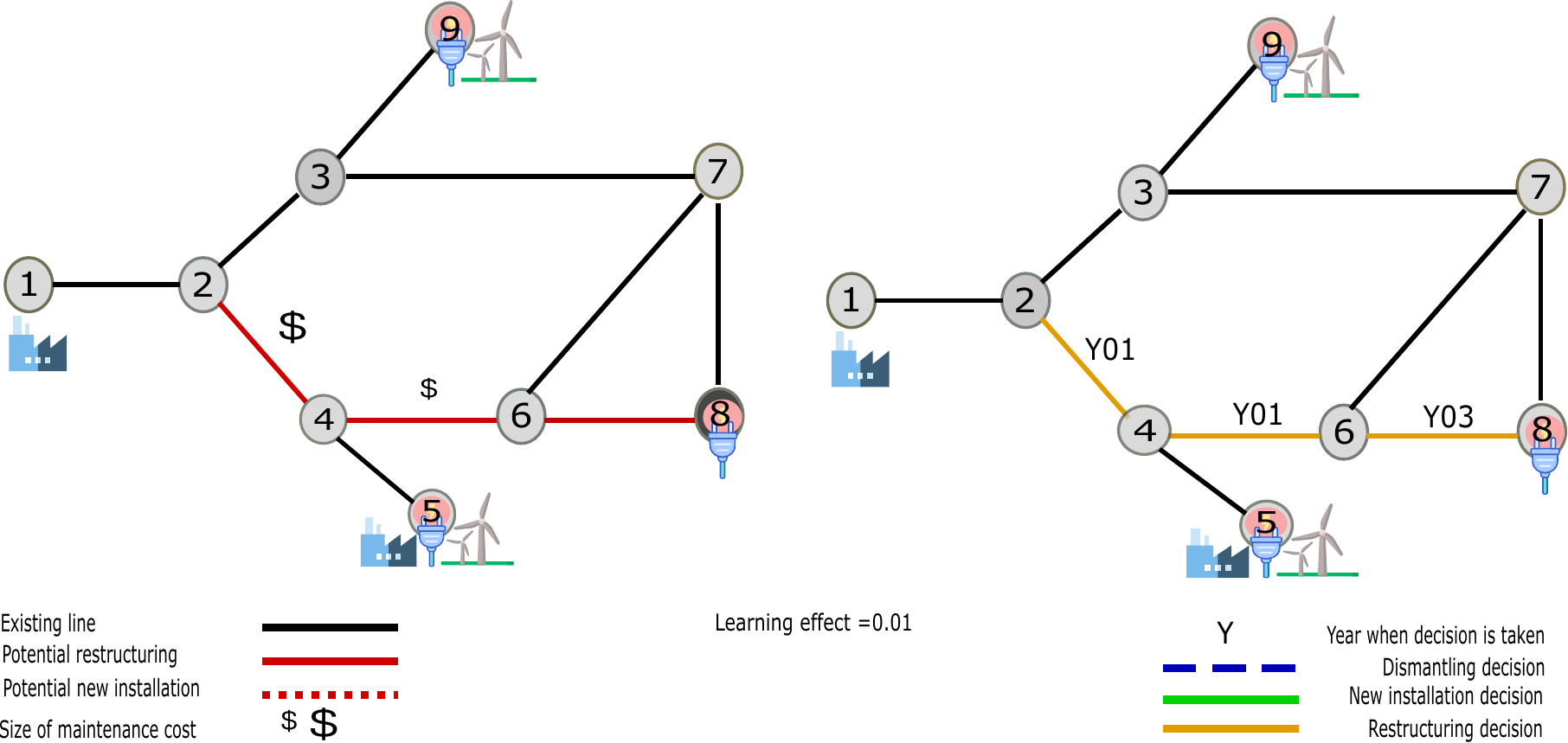}
\caption{Restructuring decision when technological learning and the maintenance cost of two cables are included}\label{Fig-Slide5}
\end{center}
\end{figure}

Let us now look at what happens when a maintenance cost is applied not only on cable 2-4 but also on cable 4-6, as shown in Figure \ref{Fig-Slide5}. Now also the cable 4-6 is dismantled on the first year, compared to Figure \ref{Fig-Slide4} where it was scheduled only on the third year.  

A question now arises: why in Figure \ref{Fig-Slide5} cables 2-4 and 4-6 are now both dismantled on the same year even though there is learning effect that should motivate the choice to schedule them in successive years? This is because the learning effect in this case is too low compared to the maintenance cost of the cables involved. Therefore it is economically more convenient to dismantle both cables together on the first year to avoid the increasing maintenance costs, compared to keep one of the cables one year more to take advantage of learning effect. 

\begin{figure}[htbp] 
\begin{center}
\includegraphics[scale=0.70]{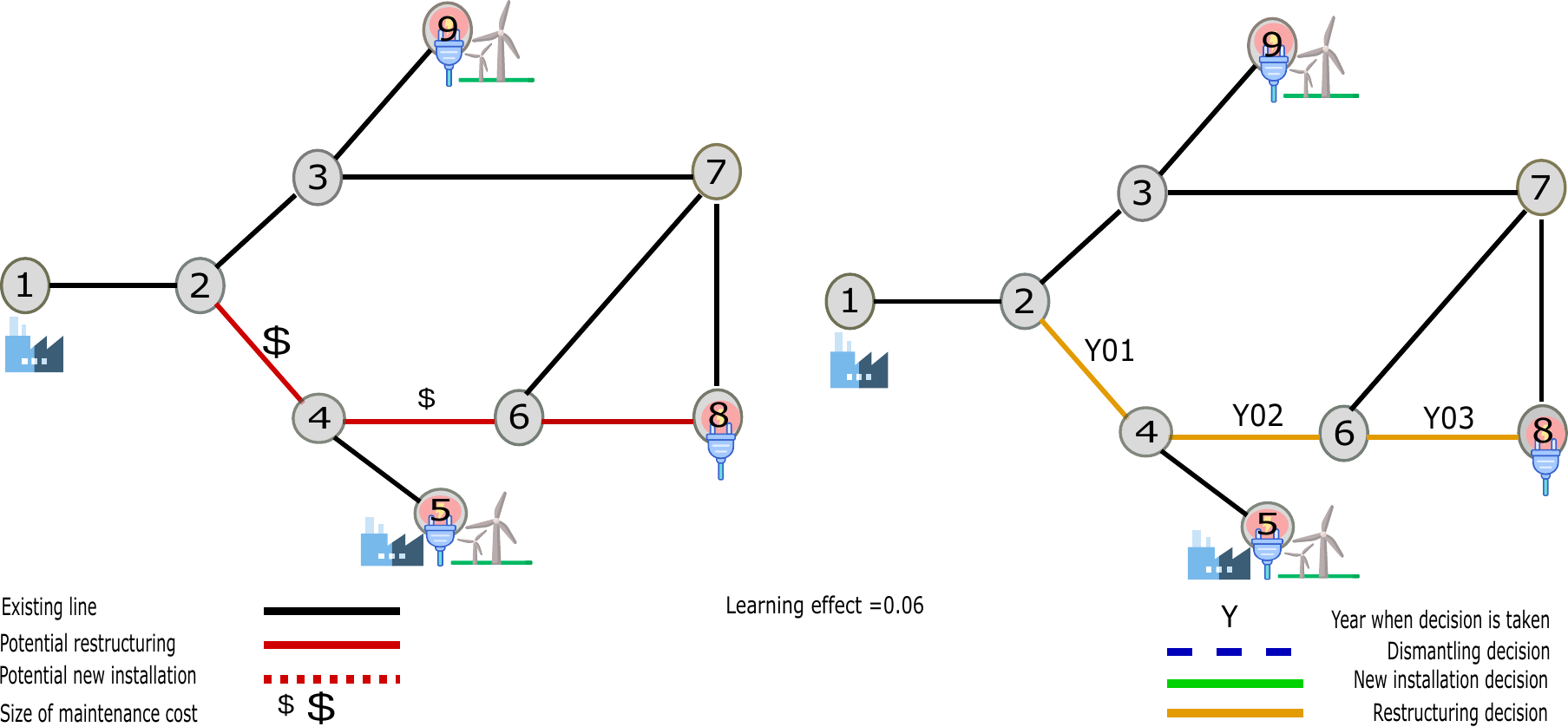}
\caption{Restructuring decision when increased technological learning and maintenance cost of two cables are included}\label{Fig-Slide6}
\end{center}
\end{figure}

Figure \ref{Fig-Slide6} shows how the decision change when a higher learning effect is applied to the case study analysed in Figure \ref{Fig-Slide5}. For this particular case, it took at least 6\% learning effect in order to motivate investment scheduling and keep one of the obsolete cables longer. In this case, the decision is to delay the replacement of cable 4-6 which indeed, is the one with the lower maintenance cost as shown in the figure.

%%%%%%%%%%%%%%%%%%%%%%%%%%%%%%%%%%%%%%%%%%%%%%%%%%%%%%%%%%%%%%%%%%%%%%%%%%%%
\subsection{Focus - Restructuring vs New Potential Installation Decisions}

The following tests aim at showing how the decisions change when also new potential installation opportunities come into the picture. We now introduce also the possibility to choose between two new different cables: one cable is bigger and more expensive (thicker green/yellow lines in the figures) while the other cable is smaller and cheaper (thinner green/yellow lines in the figures). Different ratios between the costs of the two cables, will lead to different investment choices.

\begin{figure}[htbp] 
\begin{center}
\includegraphics[scale=0.70]{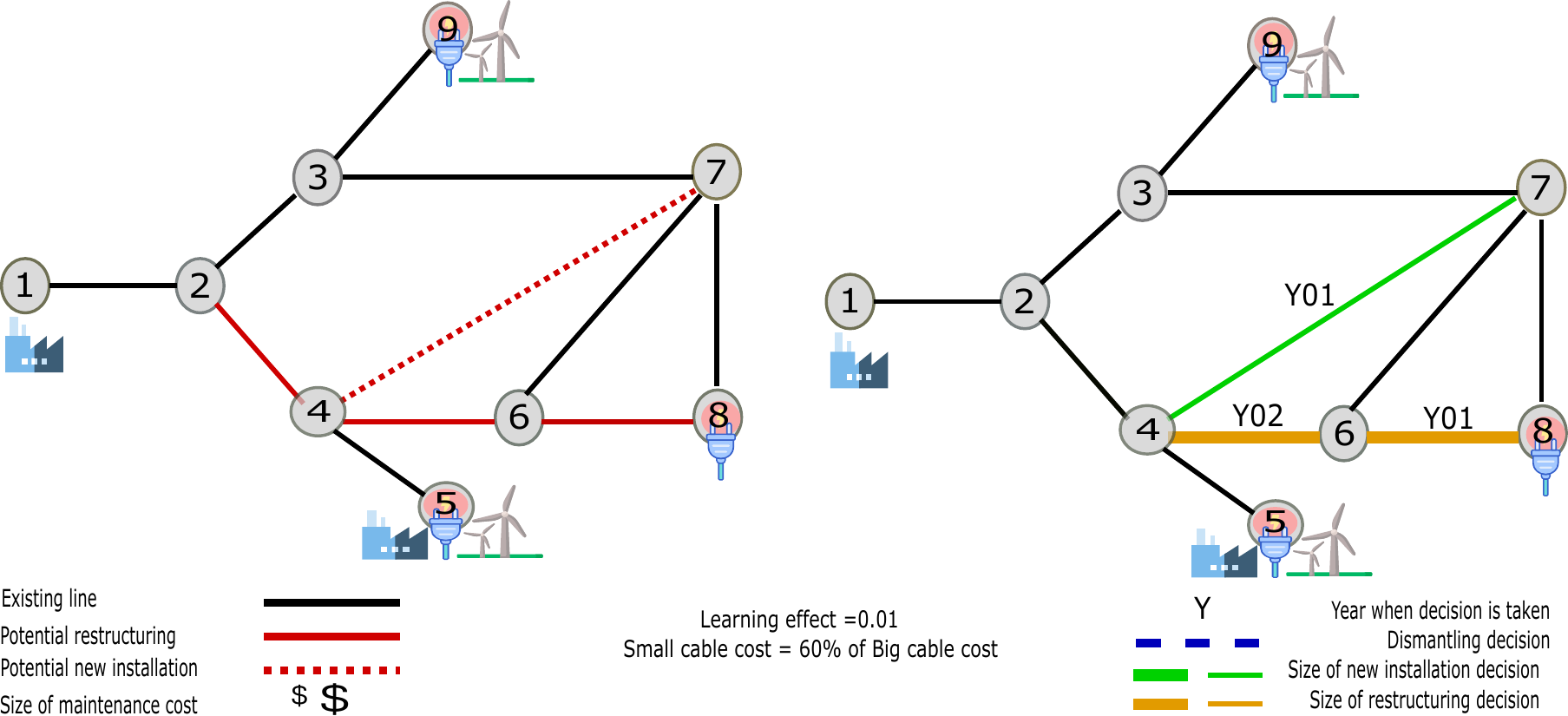}
\caption{Investment decisions when restructuring and new potential installation are included. Smaller cable cost is 60\% of the bigger cable cost}\label{Fig-Slide7}
\end{center}
\end{figure}

\begin{figure}[htbp] 
\begin{center}
\includegraphics[scale=0.70]{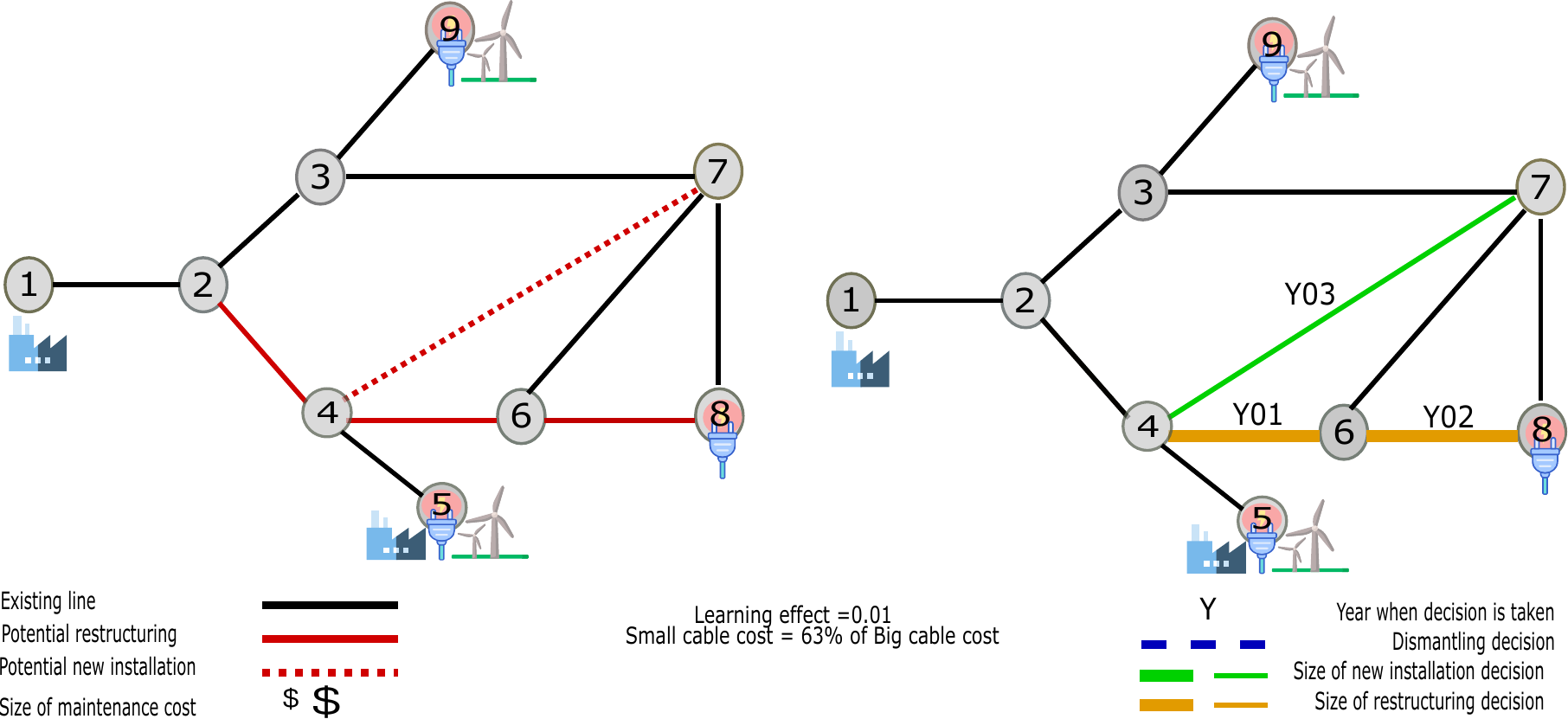}
\caption{Investment decisions when restructuring and new potential installation are included. Smaller cable cost is 63\% of the bigger cable cost}\label{Fig-Slide8}
\end{center}
\end{figure}

\begin{figure}[htbp] 
\begin{center}
\includegraphics[scale=0.70]{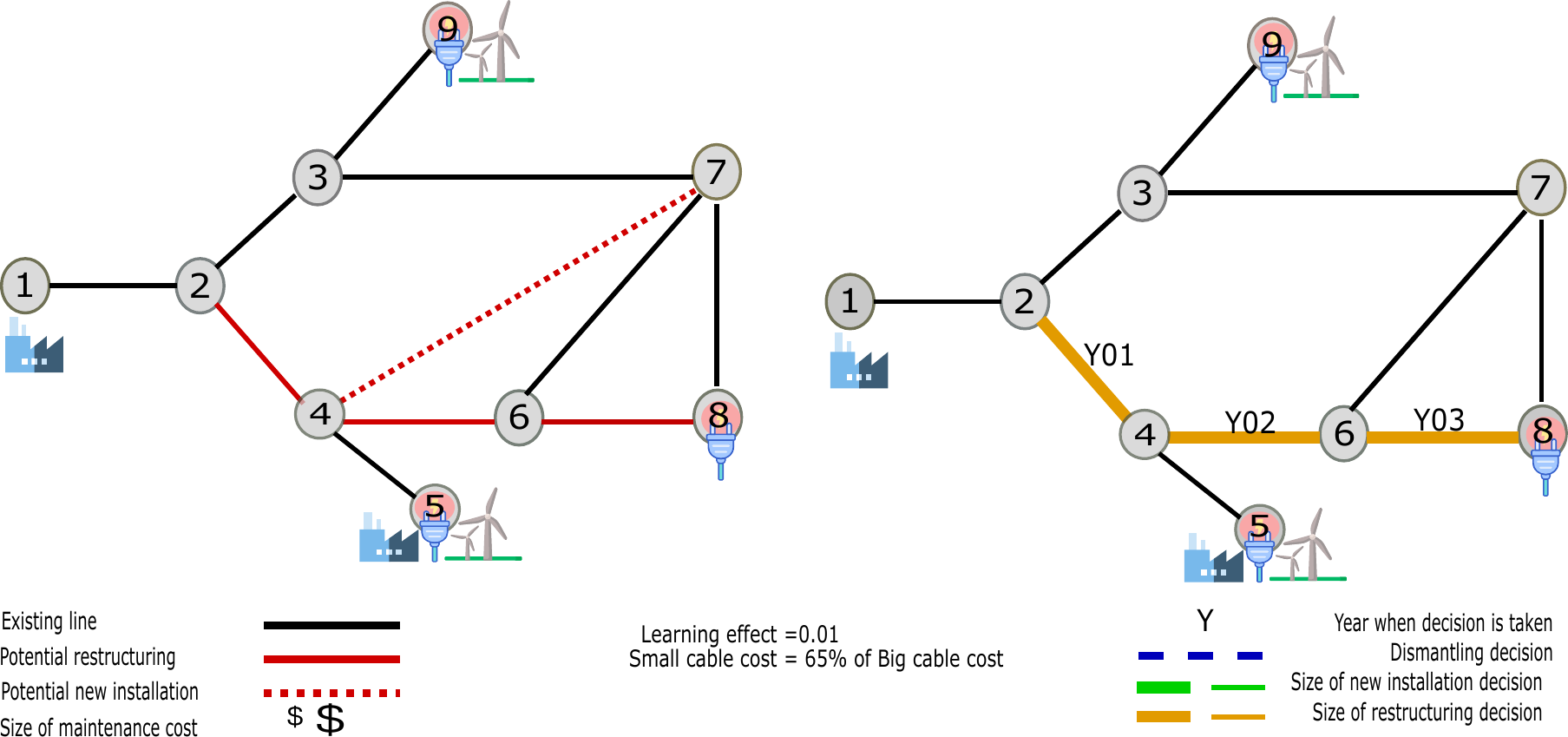}
\caption{Investment decisions when restructuring and new potential installation are included. Smaller cable cost is 65\% of the bigger cable cost}\label{Fig-Slide9}
\end{center}
\end{figure}

\begin{figure}[htbp] 
\begin{center}
\includegraphics[scale=0.70]{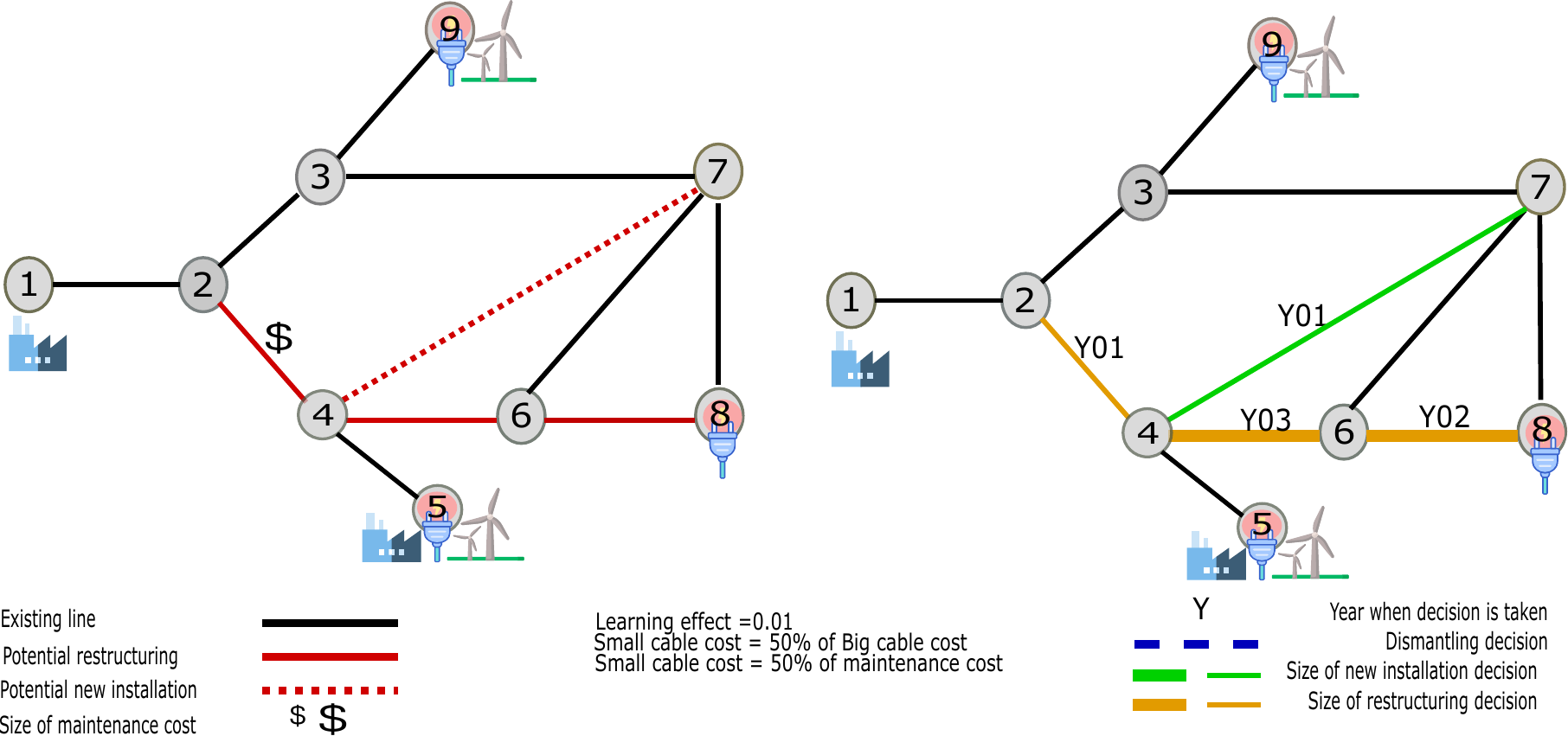}
\caption{Investment decisions when restructuring and new potential installation are included as well as the maintenance cost of a cable. Smaller cable cost is lower than 50\% of the bigger cable cost}\label{Fig-Slide10}
\end{center}
\end{figure}

Figure \ref{Fig-Slide7} should be compared with Figure \ref{Fig-Slide3}.
When the opportunity to choose between restructuring and new potential installation arises, it is possible to observe a change in the decision making process. In this particular case, compared to Figure \ref{Fig-Slide3}, cable 2-4 is no longer replaced, but a new potential cable is built between nodes 4 and 7. It is important to note that the new potential cable installed on arc 4-7 is the smaller and cheaper one. The same cable could not have been chosen for the restructuring of arc 2-4 because in that case it was necessary to reinforce the connection with a bigger and more expensive cable. The new connection 2-4 is adding new capacity to the grid. This additional new capacity, combined with the existing capacity of cable 2-4, is enough to satisfy the increasing demand, and cheaper compared to replacing the connection 2-4 with a bigger and more expensive cable.

Figures \ref{Fig-Slide8} and \ref{Fig-Slide9}, show how the investment choice changes when different ratios between the costs of the two cables are applied. As long as the smaller cable is cheaper, it is convenient to perform a new installation instead of replacing an existing cable with a more expensive one. But as the smaller cable costs increase, such choice becomes less convenient. In Figure \ref{Fig-Slide8} the smaller cable cost is assumed as 63\% of the bigger cable cost, compared to Figure \ref{Fig-Slide7} where it was set to 60\%. This affects the scheduling of the investment decisions. Due to the higher cost of the smaller cable, the new potential installation is now delayed to the 3rd year, also to take advantage of the learning coming from the restructuring of the other two obsolete cables. 

\begin{figure}[htbp] 
\begin{center}
\includegraphics[scale=0.70]{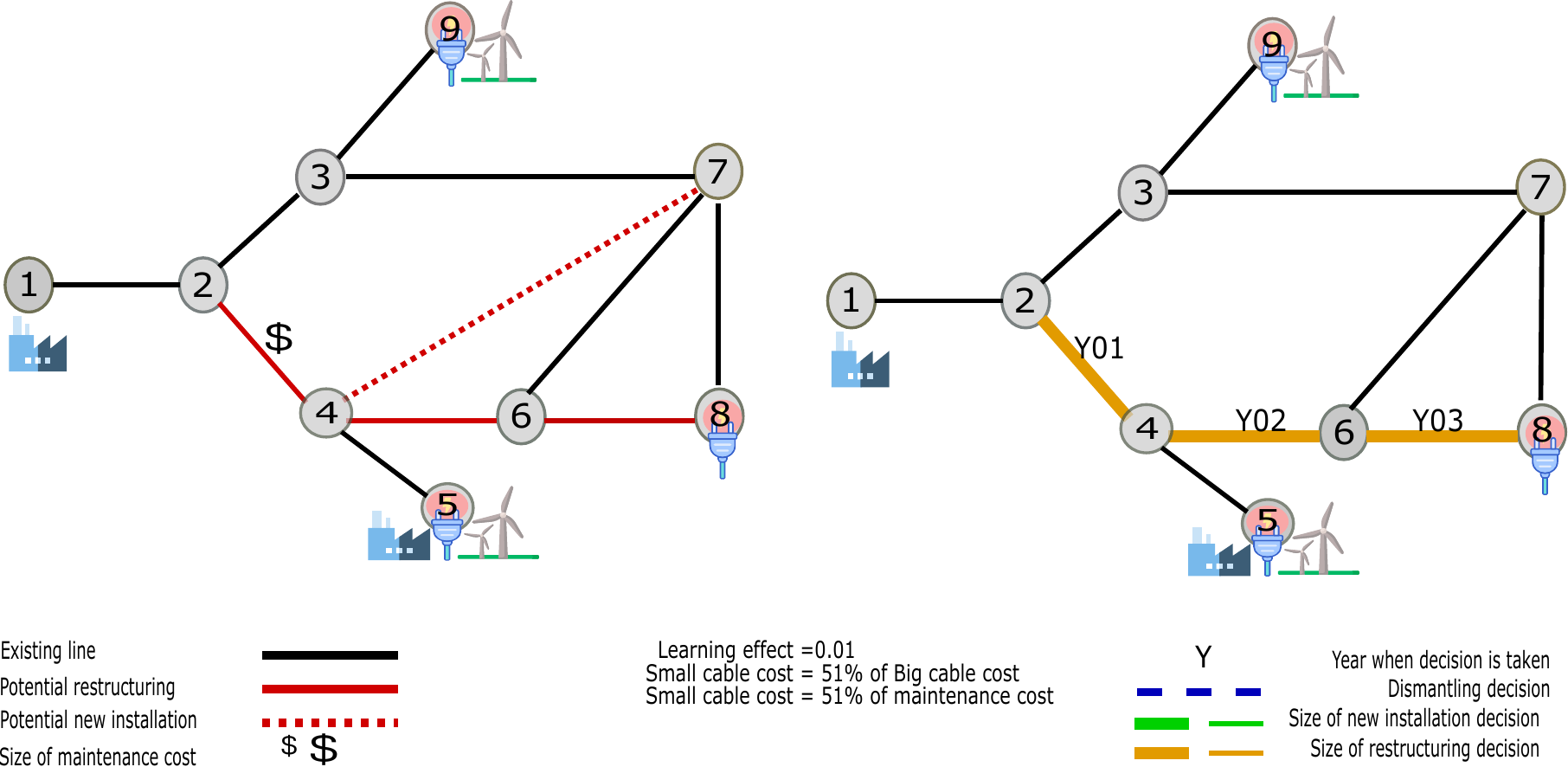}
\caption{Investment decisions when restructuring and new potential installation are included as well as the maintenance cost of a cable. Smaller cable cost is greater than 50\% of the bigger cable cost}\label{Fig-Slide11}
\end{center}
\end{figure}

In Figure \ref{Fig-Slide9} the smaller cable cost is assumed as 65\% of the bigger cable cost (hence higher than the previous tests). Now the cost of building a new connection with such a cable is too high compared to the cost of replacing the existing cable 2-4 with a new bigger one.

Figures \ref{Fig-Slide10} and \ref{Fig-Slide11} show how the investment choice changes when also maintenance cost of cable 2-4 is involved. As there is now also a maintenance cost involved in arc 2-4, the potential installation on arc 4-7 has to be done together with an additional restructuring of arc 2-4. Compared to Figure \ref{Fig-Slide7}, it is now worthy to invest on the restructuring of cable 2-4 in addition to the new connection 4-7, because cable 2-4 has a maintenance cost. However, the restructuring is performed by installing the smaller cheaper cable available.

At the same time, it was also observed that the new potential installation on arc 4-7 was performed when the smaller cable cost went down to 50\% of the bigger cable cost. For higher costs, the additional potential installation was not performed. As shown in \ref{Fig-Slide11}, when the smaller cable cost was greater than 50\% of the bigger cable cost, then only a restructuring of cable 2-4 was performed. In this case, the bigger more expensive cable was chosen for restructuring.

%%%%%%%%%%%%%%%%%%%%%%%%%%%%%%%%%%%%%%%%%%%%%%%%%%%%%%%%%%%%%%%%%%%%%%%%%%%%%%%%%%%%%%%%%
\subsection{Focus - Restructuring vs New Potential Installation vs Dismantling Decisions}

The following tests will show how the investment decisions change when more choices are available in terms of potential installation, and when also the possibility to dismantle existing obsolete cables is introduced. Different ratios between the dismantling cost and the smaller cable cost, will lead to different decisions.

\begin{figure}[htbp] 
\begin{center}
\includegraphics[scale=0.70]{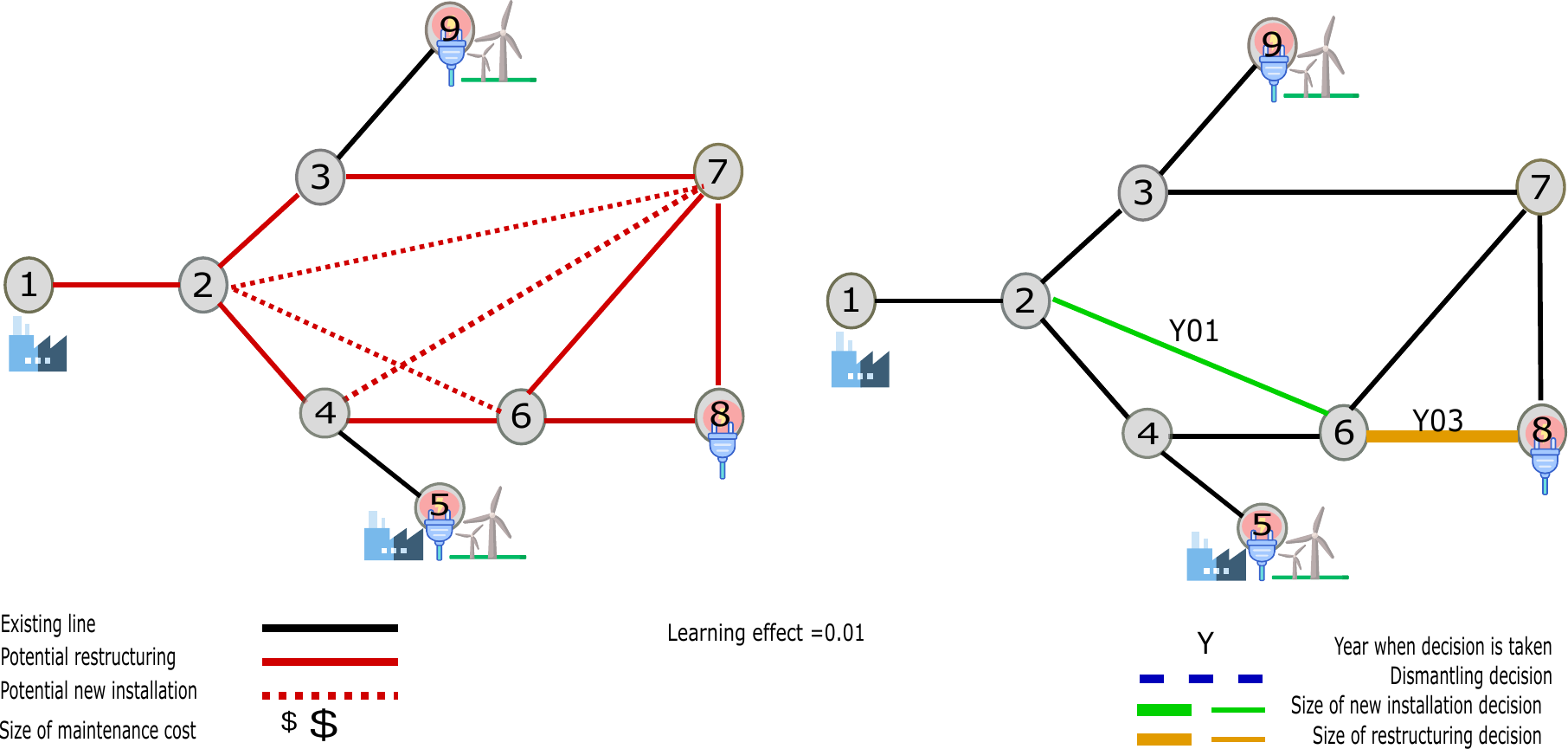}
\caption{Investment decisions when more options for new potential installation are included.}\label{Fig-Slide12}
\end{center}
\end{figure}

Figure \ref{Fig-Slide12} shows how the investment decisions change compared to Figure \ref{Fig-Slide7} when more connections for potential installations are available. Compared to Figure \ref{Fig-Slide7}, it is now possible to fulfill the future demand increment by performing just one cable restructuring and one new potential installation, (instead of two cables restructuring together with one new potential installation as it happened in Figure \ref{Fig-Slide7}). A new potential installation on arcs 2-6 instead of arcs 4-7 allows a better power flow within the network that can better fulfill the increasing demand at a cheaper cost.

\begin{figure}[htbp] 
\begin{center}
\includegraphics[scale=0.70]{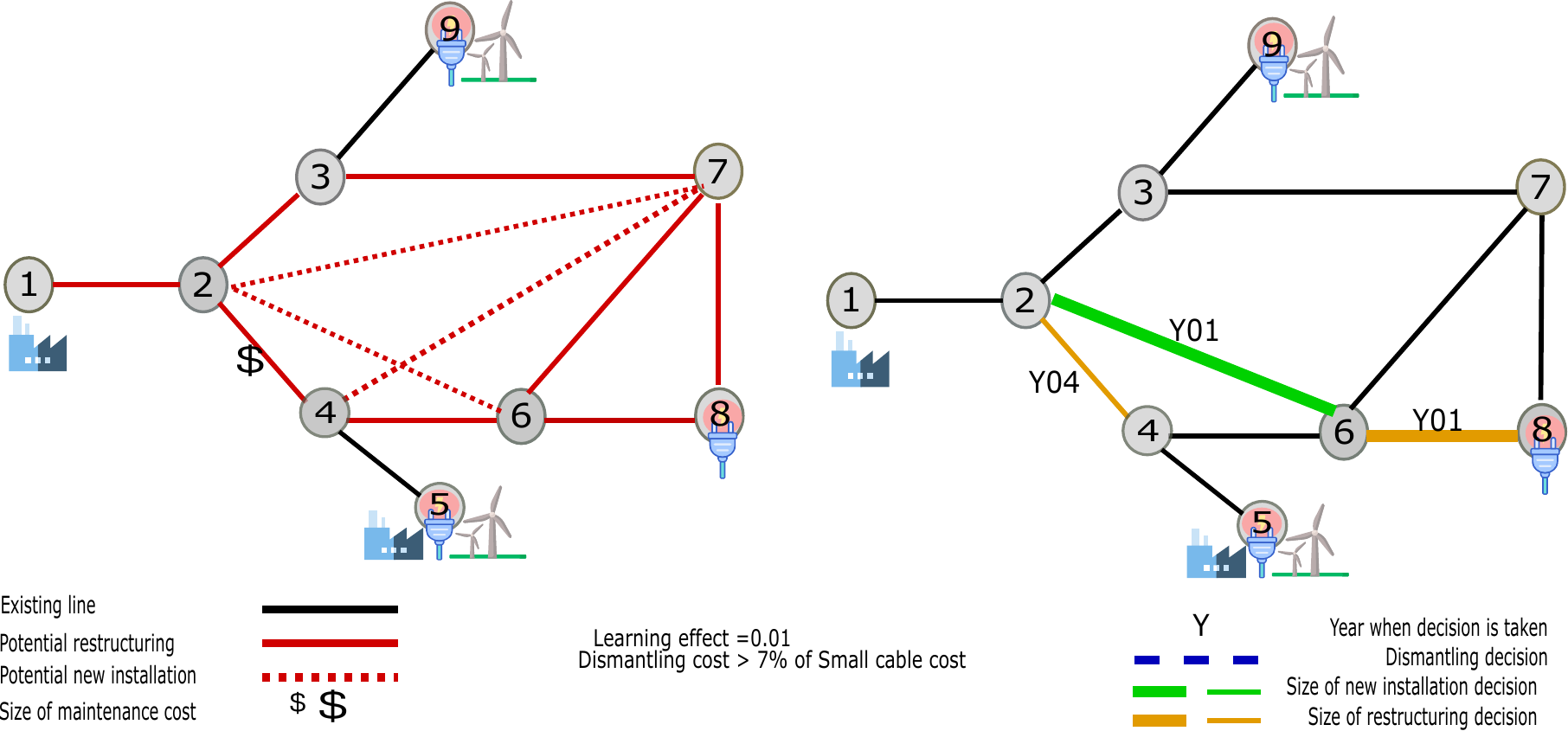}
\caption{Investment decisions when restructuring, new potential installation and dismantling possibilities are included. Maintenance cost of cables is also applied on cable 2-4. Dismantling cost is greater than 7\% of the smaller new cable.}\label{Fig-Slide13}
\end{center}
\end{figure}

\begin{figure}[htbp] 
\begin{center}
\includegraphics[scale=0.70]{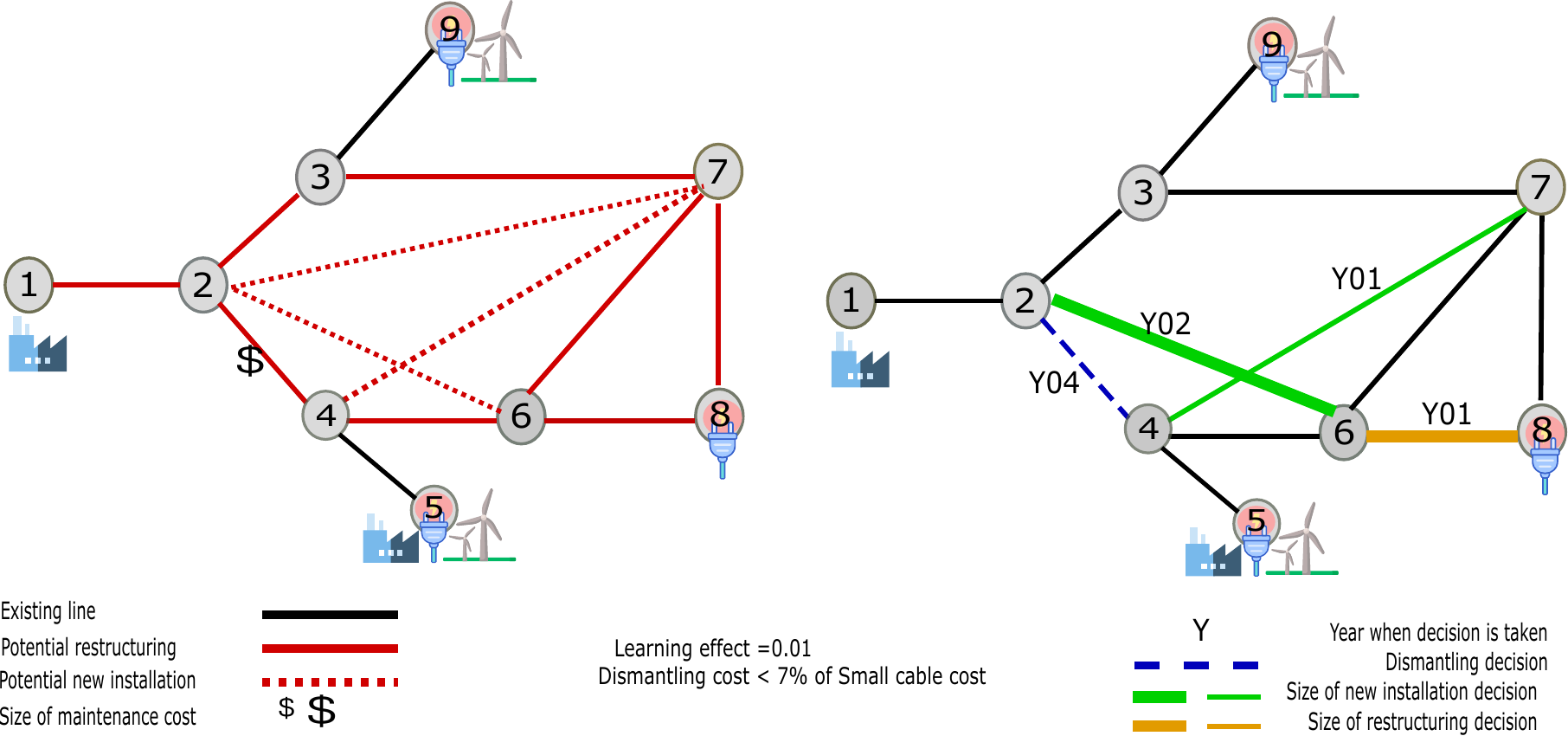}
\caption{Investment decisions when restructuring, new potential installation and dismantling possibilities are included. Maintenance cost of cables is also applied on cable 2-4. Dismantling cost is lower than 7\% of the smaller new cable.}\label{Fig-Slide14}
\end{center}
\end{figure}

Figures \ref{Fig-Slide13} and \ref{Fig-Slide14} shows how the decisions change when maintenance cost of cable 2-4 is involved as well as the possibility to dismantle existing cables. For this particular case it was found that, for high dismantling costs (greater than 7\% of the smaller cable cost), then the restructuring of cable 2-4 was preferred compared to its dismantling. While for lower dismantling costs (lower than 7\% of the smaller cable cost), then the restructuring of cable 2-4 was no longer worthy and it was more economical to dismantle it to get rid of the maintenance cost. This was leading to an additional new potential installation also on arc 4-7.

\begin{figure}[htbp] 
\begin{center}
\includegraphics[scale=0.70]{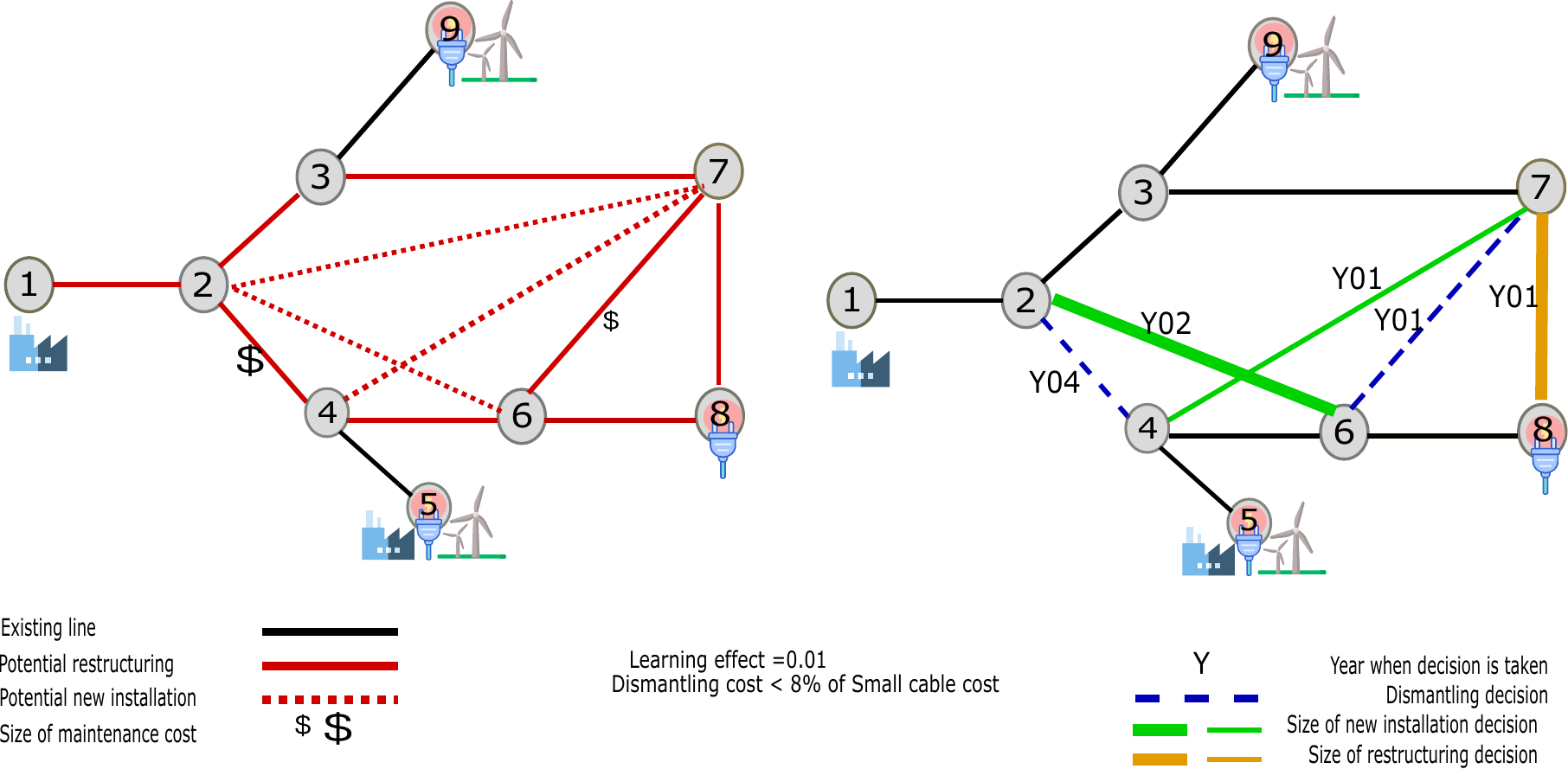}
\caption{Investment decisions when restructuring, new potential installation and dismantling possibilities are included. Maintenance cost of cables is also applied on both cable 2-4 and 6-7.}\label{Fig-Slide15}
\end{center}
\end{figure}

Figure \ref{Fig-Slide15} shows that, when a maintenance cost is applied to both cables 2-4 and 6-7, then both cables are dismantled. At the same time, potential installation on arc 2-6 and 4-7 are needed together with a restructuring of cable 7-8. 

\begin{figure}[htbp] 
\begin{center}
\includegraphics[scale=0.70]{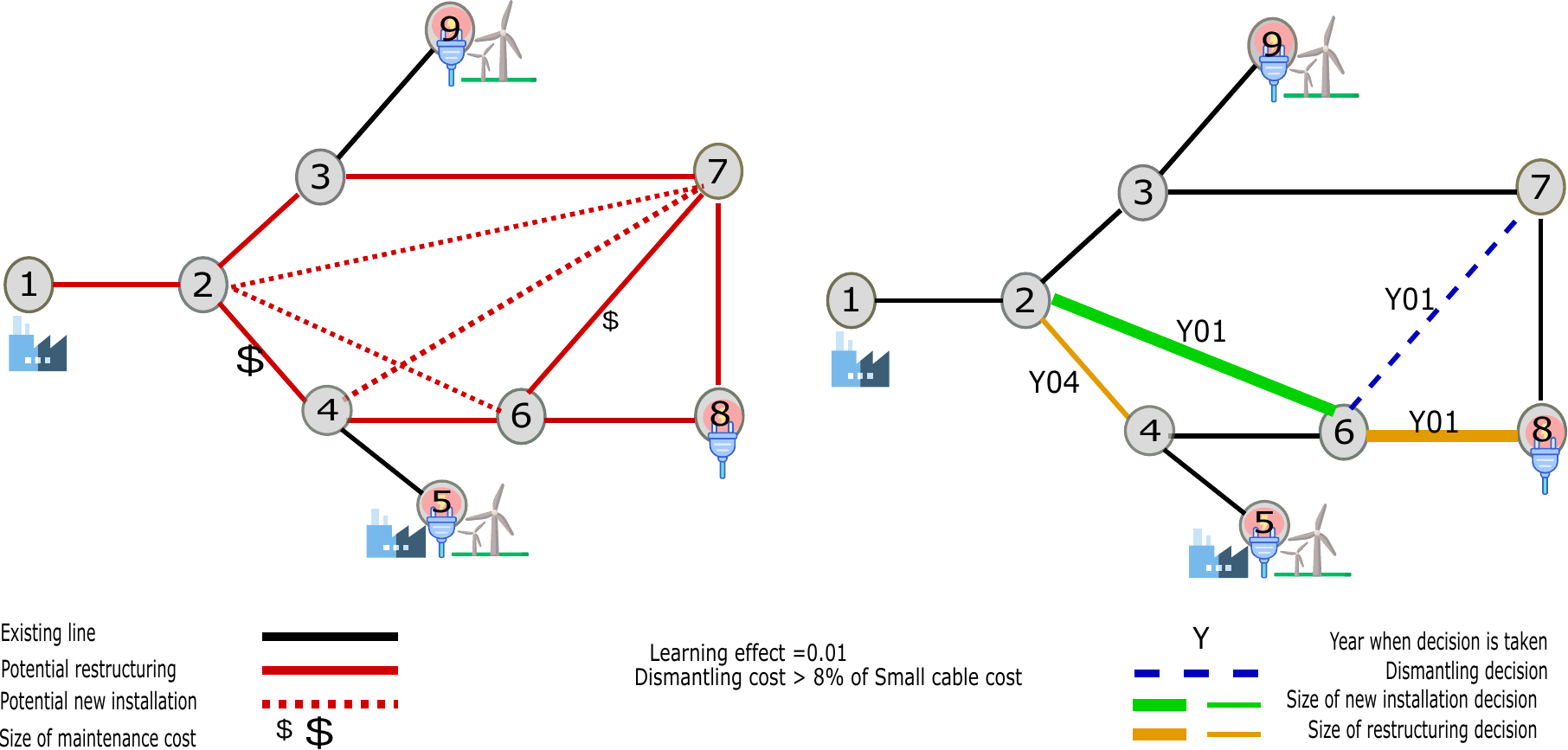}
\caption{Investment decisions when restructuring, new potential installation and dismantling possibilities are included. Maintenance cost of cables is also applied on both cable 2-4 and 6-7. }\label{Fig-Slide16}
\end{center}
\end{figure}

Figure \ref{Fig-Slide16}, compared to Figure \ref{Fig-Slide15}, shows how the decision changes when the dismantling cost is higher. Now only cable 6-7 is dismantled, while for cable 2-4 it is more economical to proceed with a restructuring by installing the smaller cheaper new cable available.

\begin{figure}[htbp] 
\begin{center}
\includegraphics[width=0.8\textwidth]{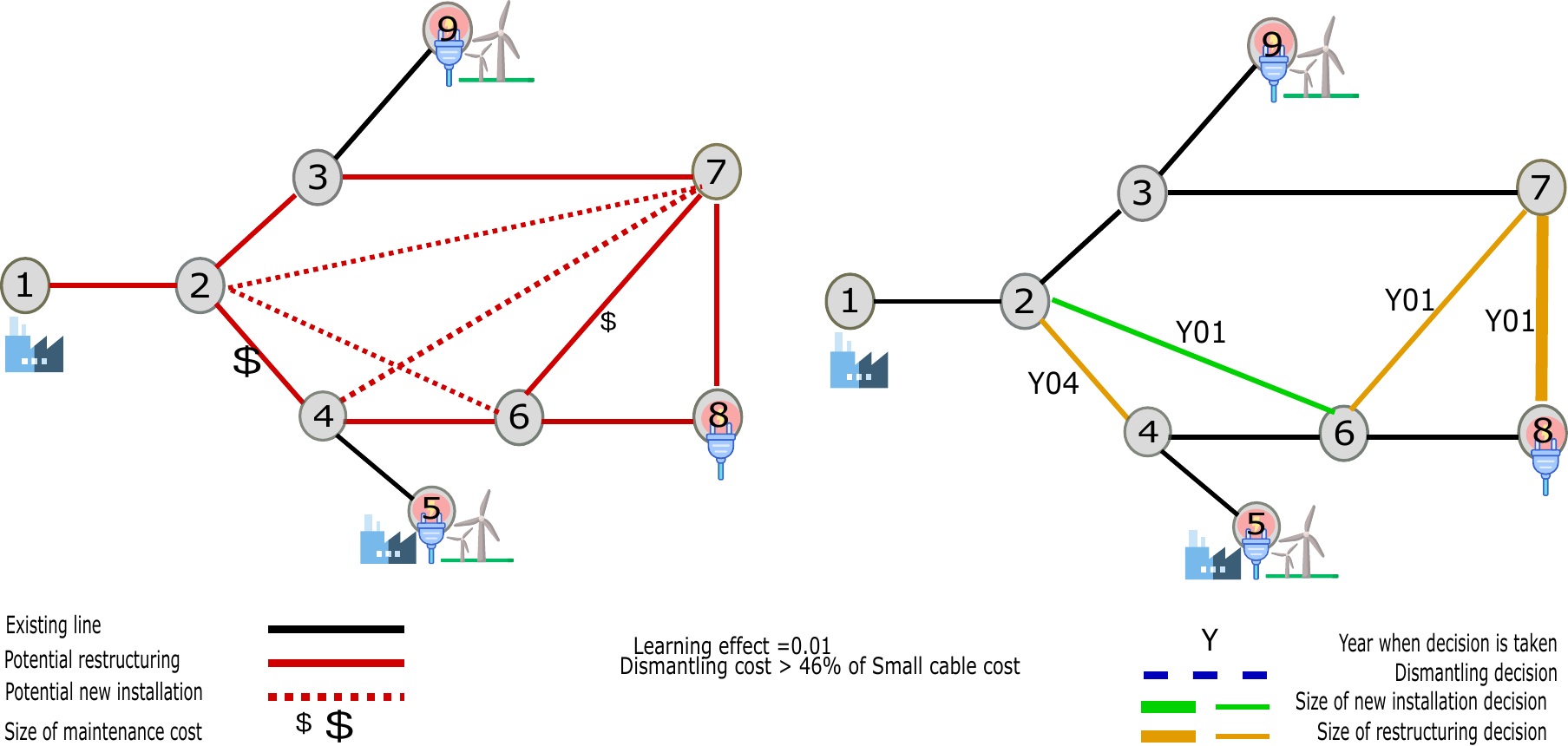}
\caption{Investment decisions when restructuring, new potential installation and dismantling possibilities are included. Maintenance cost of cables is also applied on both cable 2-4 and 6-7. Dismantling cost is assumed as very high.}\label{Fig-Slide17}
\end{center}
\end{figure}

Finally, Figure \ref{Fig-Slide17} shows the results when a very high dismantling cost is applied and therefore restructuring of obsolete cables is preferred compared to dismantling. High dismantling costs can be associated to environmental charges for disposal or particular accessibility issues to certain areas that may be remote or more difficult to handle than others.

%%%%%%%%%%%%%%%%%%%%%%%%%%%%%%%%%%%%%%%%%%%%%%%%%%%%%%%%%%%%%%%%%%%%%%%%%%%%%%%%%%%%%%%%
\section{Conclusions and Research Directions}\label{Conclusions}
In this paper we have discussed a multihorizon approach for the reliability oriented network restructuring problem by including additional key modelling features: technological learning, construction time and cables maintenance costs, as well as the possibility to consider future demand increment.

The computational experiments and analyses show that, such features, have key impact on the final decision making process. They not only affect the restructuring decisions here and now, but also the investment scheduling in the forthcoming years.

The ability to put together historical data to forecast the future maintenance costs of existing cables and apparatus is therefore of high importance when planning investments in network expansion, reinforcement and reconfiguration. The choice between keeping an existing connection as it is or changing it, is highly dependent on the forecast behaviour of the existing apparatus. Hence maintenance cost functions for the existing electrical apparatus have to be properly defined, forecast and calculated. Future research directions should be focused on a proper definition of cables forecast maintenance costs through machine learning. A proper definition of such maintenance costs would provide finer input dataset to the optimisation models and therefore much better decision making processes.

Technological learning has also turned out very important for this kind of problems, especially when it comes to schedule investments in restructuring and reconfiguration of different corridors of the grid. However, the learning rates of the technologies are uncertain. Moreover, this is the first time that technological learning is introduced within an optimisation model for the reliability oriented network restructuring problem. Therefore, the range of a learning rate of learning-by-doing for this particular type of problems should be further investigated. In order to capture this
aspect, a stochastic programming approach can be also applied. With uncertain learning rates, a more sophisticated path of cables installations might be followed, and different scheduling choices might take place. 
The ability to define proper technological learning coefficients for different activities within the reliability oriented network restructuring problem, can improve the investment scheduling. 

%\section*{Acknowledgement}
%This work was partly supported by the Estonian Research Council grant PUTJD915.

%\section*{References}
\bibliographystyle{unsrt}  % (uses file "plain.bst")amsplain
\bibliography{references}

\end{document}